\documentstyle[amssymb,amsfonts]{amsart}

\def\R{{\hbox{\bf R}}}
\def\C{{\hbox{\bf C}}}

\font \roman = cmr10 at 10 true pt

\renewcommand{\H}{\mbox{{\bf H}}}

\def\allt#1{%
\smash{
\vtop{%
     \ialign{%
        ##\crcr
        $\hfil\displaystyle{\tilde \forall}\hfil$\crcr%
        \noalign{\kern1.5pt\nointerlineskip}
        $\hfil\!\!#1\hfil$\crcr\noalign{\kern1.5pt}
        }
       }
      } \hbox{$\vphantom{#1}$}
     }

\def\be#1{ \begin{equation}\label{#1} }

\def\bas{\begin{align*}}
\def\eas{\end{align*}}
\def\bi{\begin{itemize}}
\def\ei{\end{itemize}}

\def\diam{{\hbox{\roman diam}}}

\def\M{{\hbox{\roman M}}}

\def\dist{{\hbox{\roman dist}}}

\def\Z{{\hbox{\bf Z}}}
\def\U{{\hbox{\bf U}}}

\def\eps{\varepsilon}
\newenvironment{proof}{\noindent {\bf Proof} }{\endprf\par}
\def \endprf{\hfill  {\vrule height6pt width6pt depth0pt}\medskip}
\def\emph#1{{\it #1}}
\def\textbf#1{{\bf #1}}
\def\g{{\frak g}}

\theoremstyle{plain}
  \newtheorem{theorem}[subsection]{Theorem}

  \newtheorem{proposition}[subsection]{Proposition}
  
  \newtheorem{lemma}[subsection]{Lemma}
  \newtheorem{corollary}[subsection]{Corollary}

\theoremstyle{remark}

\theoremstyle{definition}
  \newtheorem{definition}[subsection]{Definition}

\include{psfig}

\begin{document}

\title[Singularity removal for Yang-Mills]{A singularity removal theorem for 
Yang-Mills fields in higher dimensions}

\author{Terence Tao}
\address{Department of Mathematics, UCLA, Los Angeles CA 90095-1555}
\email{tao@@math.ucla.edu}
\thanks{TT is a Clay Prize Fellow and is supported by a grant from the Packard Foundation.}

\author{Gang Tian}
\address{Department of Mathematics, Massachusetts Institute of Technology,
Cambridge, MA 02139}
\email{tian@@math.mit.edu}
\thanks{GT is supported by a NSF grant and a Simons fund.}

\begin{abstract}  In four and higher dimensions, we show that any stationary admissible
Yang-Mills field can be gauge transformed to a smooth field if the $L^2$ norm
of the curvature is sufficiently small.  
There are three main ingredients. The first is Price's monotonicity formula, which
allows us to assert that the curvature is small not only in the $L^2$ norm, but
also in the Morrey norm $M_2^{n/2}$.  The second ingredient is a new inductive
(averaged radial) gauge construction and
truncation argument which allows us to approximate a singular gauge as a weak
limit of smooth gauges with curvature small in the Morrey norm.
The second ingredient is variant of Uhlenbeck's lemma, allowing one to place a
smooth connection into the Coulomb gauge whenever the Morrey norm of the
curvature is small; This variant was also proved
independently by Meyer and Riviere \cite{riviere}. It follows easily from this variant that 
a $W^{1,2}$-connection can be placed in the Coulomb gauge if it can be approximated by smooth connections
whose curvatures have small Morrey norm.
\end{abstract}

\maketitle

\section{Introduction}\label{intro-sec}

The purpose of this paper is to investigate the small-energy behavior of weakly
Yang-Mills fields in $\R^n$ for
$n \geq 4$, and in particular to extend the singularity removal theorem of
Uhlenbeck \cite{uhlenbeck-ym} to higher dimensions.

Fix $n \geq 4$, and let $\Omega$ be some bounded domain in $\R^n$; typically we
shall restrict our attention to the
cubes $\Omega = [-1,2]^n$ or $\Omega = [0,1]^n$.

Let $G$ be a fixed finite-dimensional compact Lie group; it will be convenient
to consider $G$ as embedded in some large
unitary group $U(N)$. Let $\g$ be the Lie
algebra of $G$.  We define a \emph{connection} on $\Omega$ to be a section $A$
of $T^* \Omega \otimes \g$ (i.e. a
$\g$-valued 1-form) which is locally $L^2$.
For any connection $A$ let
\be{f-def}
F(A) := dA + A \wedge A
\end{equation}
denote the curvature of $A$.  Since $A$ is locally $L^2$, $F(A)$ makes sense as
a ($\g$-valued 2-form) distribution.

A \emph{gauge transformation} is a sufficiently regular\footnote{As a bare
minimum, one should have $\sigma$ in the
Sobolev space $W^{1,2}_{loc}$ and also in $L^\infty$, in order for $\sigma(A)$
to be locally in $L^2$.  In practice
we shall have significantly more regularity than this.} map $\sigma: \Omega \to
G$.  This group acts on connections by the formula
\be{e-def}
\sigma(A) := \sigma \cdot A \cdot \sigma^{-1} - d\sigma \cdot \sigma^{-1}.
\end{equation}
We call $A$ and $\sigma(A)$ \emph{gauge equivalent}.

Let $A$ be a connection on the cube $[-1,2]^n$.  We say that $A$ is a
\emph{smooth Yang-Mills connection} on
$[-1,2]^n$ if $A$ is smooth and solves the PDE
\be{ym}
 d_* F(A) -\ast [A, \ast F(A)] = 0,
\end{equation}
where $\ast$ denotes the Hodge operator and $d_* := \ast d \ast$;
one may verify that this condition is invariant under gauge transformations.
Formally, Yang-Mills connections are critical points of the energy
functional $\int |F(A)|^2$.  Following \cite{tian}, we say that
$A$ is an \emph{admissible
Yang-Mills connection} if it is a smooth Yang-Mills connection
outside a closed subset
$S\subset [-1,2]^n$ of finite
$(n-4)$-dimensional Hausdorff measure and $\int |F(A)|^2 < \infty$.
We will call $S$ the \emph{singular set} of $S$.

It follows from analysis in \cite{tian} that weak limits of smooth Ynag-Mills connections with curvature uniformly
$L^2$-bounded are admissible Yang-Mills connections.

Following \cite{tian}, we call an admissible Yang-Mills connection $A$ \emph{stationary} if for any vector
field $X=X^i {\partial
\over \partial x_i}$ with compact support in $(-1,2)^n$, we have
$$
\int_{[-1,2]^n} \left( |F_A|^2 {\rm div}(X) - 4 F_{\alpha \beta}F_{\alpha
\gamma}
{\partial X^\beta\over\partial x_\gamma}\right ) = 0,
$$
where $F_A = F_{\alpha \beta} dx_\alpha\wedge dx_\beta$ and we use the usual
summation conventions. It follows from a monotonicity formula of Price
\cite{price} that
$r^{4-n} \int_{B(x,r)}|F_A|^2$ is monotone non-decreasing for any stationary
Yang-Mills connection $A$
(cf. \cite{tian}).

The main result of this paper is the following singularity removal theorem
for small energy stationary admissible Yang-Mills connections:

\begin{theorem}\label{abstract3}  Let $A$ be a stationary admissible Yang-Mills
connection on $[-1,2]^n$ with singular set $S$
which obeys the smallness condition
\be{F-small}
\int_{[-1,2]^n} |F(A)|^2 \leq \eps.
\end{equation}
Then, if $0 < \eps \ll 1$ is sufficiently small (depending only on $n, G, N$),
there is a gauge transformation $\sigma$ on $[0,
1]^n \backslash S$ such that $\sigma(A)$ extends to a smooth connection over
all of $[0,1]^n$.  Indeed, we have the
uniform bounds
\be{ns}
| \nabla^j \sigma(A)(x)| \leq C_j \eps
\end{equation}
for\footnote{Here and in the sequel, $C$ denotes various absolute constants
depending only on $n$, $G$, and $N$.} all $x
\in [0,1]^n$ and $j = 0, 1, \ldots$.
\end{theorem}

When $n=4$, the stationary property is automatic for admissible Yang-Mills
connections.
Hence, the above theorem generalizes the results of Uhlenbeck
(\cite{uhlenbeck},
\cite{uhlenbeck-ym}), who proved the above removable singularity theorem
in four dimensions. A proof of this theorem was given in \cite{tian} under the assumption of existence of a 
good gauge. The construction of this good gauge was unknown then even in the case that $S$ is
a smooth submanifold \footnote{This ws pointed out to the second author by K. Uhlenbeck}.
The main technical part of this paper is to fill in this gap by constructing
a Coulomb gauge for any stationary admissible Yang-Mills connection with small
$L^2$-norm of curvature. This turns out to be highly nontrivial.

The cubes $[-1,2]^n$ and $[0,1]^n$ can of course be rescaled, however one
should caution that in the higher-dimensional
case $n > 4$, the energy $\int |F(A)|^2$ is not invariant under scaling, so if
one were for instance to replace $[-1,2]^n$
by $[-r, 2r]^n$ in \eqref{F-small} then the right-hand side should be replaced
by $\eps r^{n-4}$.  It will also
be clear from the proof that the underlying space $\R^n$ can be replaced by a
smooth $n$-dimensional manifold.

Let $A$ be any stationary admissible Yang-Mills connection.
By the monotonicity formula of Price \cite{price}, the integral
$r^{4-n}\int_{B(x,r)}|F_A|^2$ is non-decreasing, so the density function
\be{density}
\Theta(A,x)= \lim_{r\to 0} r^{4-n}\int_{B(x,r)} |F_A|^2
\end{equation}
exists for any $x\in \Omega$.
The above theorem implies that $x$ is a singularity of
$A$ modulo all gauge transformations if and only if
$\Theta(A,x)\ge \epsilon$.
In fact, one can give a better lower bound for the density at a genuine
singular point $x$:
\be{thetabound}
\Theta(A,x) \ge \min\{ \inf_B c_{n,4}
\int_{S^4}|F_B|^2,\inf_{5\le k\le n}\inf_{B'}
\frac{c_{n,k}}{k-4} \int_{S^{k-1}} |F_{B'}|^2 \},
\end{equation}
where $B$ and $B'$ run over all non-flat Yang-Mills connections on $S^4$ and
$S^{k-1}$, respectively, and furthermore,
\be{cnk}
c_{n,k}=
\int_{B^{n-k}(0,1)} (\sqrt{1-r^2})^{k-4} dv,
\end{equation}
where $B^{n-k}(0,1)$ denotes the unit ball in $\R^{n-k}$.
Its proof can be outlined
as follows: If $\Theta(A,x)$ is smaller than the given number,
then it follows from results in \cite{tian} that modulo gauge transformations,
by taking a subsequence if necessary, scaled connections $\lambda A(x+\lambda
(y-x))$ converge to a Yang-Mills connection $A_\infty$ on $\R^n\backslash
S$ which is simply the homogeneous extension of a Yang-Mills connection
on $S^{n-1}$, where $S$ is a closed subset with $n-4$-dimensional
Hausdorff measure zero and which is invariant under scalings.
Then the claim follows from direct computations and induction
on dimensions.

The proof of Theorem \ref{abstract3} is somewhat lengthy and proceeds in
several stages, which we now
describe.

The first step is to recall that stationary
admissible Yang-Mills connections are smooth outside of a small set 
(a compact singular set $S$ of codimension at least
four).  Furthermore, thanks to Price's monotonicity
formula \cite{price}, the curvature $F(A)$ is not only small in the $L^2$ norm
\eqref{F-small}, but also small
in a certain \emph{Morrey space} $M^{n/2}_2([0,1]^n)$, defined below.  This
will be important because in the higher dimensional
case $n > 4$ the $L^2$ norm is not scale invariant, but the Morrey norm is.
It follows from this step that $\rho^2 F(A)$ is uniformly bounded outside
$S$, in fact, it is small near $S$, where $\rho$ denotes the distance from the set $S$.

Next we show that smooth connections with small $L^2$-norm of curvature can be placed 
in the \emph{Coulomb gauge} $d_* A = 0$, following the approach of Uhlenbeck \cite{uhlenbeck}.  
This can be done because these smooth connections have small curvature in the Morrey norm.  
To do this we need to generalize Uhlenbeck's lemma on
Coulomb gauges from Lebesgue spaces to Morrey spaces, which turns out to be
relatively standard. This generalization has also been achieved independently by Riviere and Meyer
\cite{riviere}). It follows easily from this step that if a $W^{1,2}$-connection
$A$ is weakly a limit of smooth connections with small $L^2$-norm of curvature, then $A$ has
a Coulomb gauge. 

In the next step, we excise the singular set $S$ by approximating the connection $A$ as a
weak limit of smooth connections. The
difficulty here is in ensuring that the smooth connections still have small
curvature. If $S$ is a union of disjoint smooth submanifolds, one can construct
such approximations by first eliminating the component of $A$ in $\rho$-direction and using 
the curvature estimate in the first step \footnote{We have obtained this as well as last step
long before we could do the next.}. However, this problem turns out to be
surprisingly non-trivial if we do not have a prior knowledge on smooth structure of $S$.
We will proceed by first performing an inductive gauge
transform, averaging various
radial gauges together, to transform the connection $A$ to one which obeys good
bounds away from the singular set $S$
(roughly speaking, we need a connection which blows up like $\eps / \dist(x,S)$).  
We then truncate this transformed connection by a cutoff function to obtain the approximating
connections.

These last two steps are the main technical parts of this paper. 

Finally, by taking limits, we can conclude from the above two that
the original Yang-Mills connection $A$ can be placed in the
Coulomb gauge. At this point one can use the Yang-Mills equation and some standard elliptic theory to obtain the
desired regularity of $A$, even across the singular set $S$. This was already done in \cite{tian}
in a different way.

\section{Notation}

In this section we lay out some notation, especially relating to the Lie group $G$, and the relationships between
connections $A$, curvatures $F$, and gauge transforms $\sigma$.  A useful heuristic\footnote{Related to this heuristic
is the following dimensional analysis: if we give distances in $\R^n$ the scaling of $length$, then $\sigma$ has units of $length^0$, $A$ has units of $length^{-1}$, and $F$ has units of $length^{-2}$, since each derivative in space effectively has the units of $length^{-1}$.  Meanwhile integration on $n$-dimensional sets (e.g. balls $B(x,r)$ or the cube $[0,1]^n$) effectively has units of $length^n$, and constants such as $C$ and $\eps$ are dimensionless.  The reader may then check that all of the estimates in this paper are dimensionally consistent.} to keep in mind is that
the curvature acts like one derivative of the connection, which in turn acts like one derivative of the
gauge transform (cf. \eqref{e-def}, \eqref{f-def}).

In this paper we use $C$ to denote various constants which depend only on the ambient dimension $n$, the Lie group $G$,
and the dimension $N$ of the unitary group $U(N)$ containing $G$.  We use $A \sim B$ to denote the estimate
$C^{-1} A \leq B \leq CA$.

We use $B(x,r) := \{ y \in \R^n: |y-x| < r \}$ to denote the open ball of radius $r$ centered at $x$.  If $E$
is a set in $\R^n$, we use $|E|$ to denote the Lebesgue measure of $E$, thus for instance $|B(x,r)| = C r^n$.
Also note that if $0 < r \leq 1$ and $x \in [0,1]^n$ then $|B(x,r) \cap [0,1]^n| \sim r^n$.

Recall that the Lie group $G$ is embedded in a unitary group $U(N)$, so that $\g$ is embedded in the vector space $u(N)$.  
In particular we have $|\sigma \cdot A \cdot \sigma^{-1}| = |A|$ for all $\sigma \in G$ and $A \in \g$, where $|A|$ denotes 
the operator norm in $u(N)$.   We use $1_G$ to denote the identity element in $G$.

If $A$ is a connection, we define $|A(x)| := (\sum_\alpha |A_\alpha(x)|^2)^{1/2}$
and $|F(A(x))| := (\sum_{\alpha,\beta} |F_{\alpha \beta}(A)(x)|^2)^{1/2}$.  
From the identity
\be{gauge-identity}
F(\sigma(A)) = \sigma F(A) \sigma^{-1}
\end{equation}
we observe that the magnitude $|F(A)|$ of the 
curvature is gauge invariant:
\be{whip}
|F(\sigma(A))| = |F(A)|.
\end{equation}
For future reference, we also record the composition law
\be{composition}
\sigma_1(\sigma_2(A)) = (\sigma_1 \sigma_2)(A).
\end{equation}

\section{Some preliminaries on Morrey spaces}

In this section we set up some basic notation, in particular the notation for Morrey spaces, and develop
some of the basic functional theory for these spaces such as fractional integration, Sobolev embedding, etc.

Suppose $A$ is a stationary 
Yang-Mills connection obeying the curvature smallness condition \eqref{F-small}.
By using Price's monotonicity formula \cite{price} as in \cite{tian}, we can improve \eqref{F-small} to the scale-invariant bounds
$$ \int_{B(x,r) \cap [0,1]^n} |F(A)|^2 \leq C \eps r^{n-4}$$
for all balls $B(x,r)$.  This is equivalent to a Morrey space estimate on $F(A)$, and motivates
introducing the following (standard) notation.

We follow the notation of \cite{taylor}:

\begin{definition}\label{morrey-def}
If $\Omega$ is a domain and $1 \leq q \leq p$, we define the Morrey spaces 
$M^p_q(\Omega)$ to be those locally $L^q$ functions (possibly vector-valued) 
whose norm
$$ \| f\|_{M^p_q} := \sup_{x_0 \in \R^n; 0 < r \leq 1} r^{n(\frac{1}{p} - 
\frac{1}{q})} (\int_{B(x_0,r) \cap \Omega} |f|^q)^{1/q}$$
is finite.  
We also define Morrey-Sobolev spaces $M^p_{q,k}$ for integers $k \geq 0$ by the 
formula
$$ \| f \|_{M^p_{q,k}(\Omega)} := \sum_{j=0}^k \| \nabla^j f 
\|_{M^p_q(\Omega)}.$$
In practice $k$ shall always be 0, 1, or 2.
\end{definition}

Thus Price's monotonicity formula gives $M^{n/2}_2([0,1]^n)$ control on $F(A)$.

The norm $M^p_q$ has the scaling of $L^p$, but the functions are only $L^q$ integrable.  
From H\"older's inequality we see that all $L^p$ functions are in $M^p_q$, but not 
conversely.  Note that the $M^p_q$ norm depends only on the magnitude of $f$.  
In particular, we see from \eqref{whip} that
\be{gauge-eq}
\| F(\sigma(A)) \|_{M^p_q} = \| F(A) \|_{M^p_q}.
\end{equation}
This gauge invariance of the Morrey norms for curvatures will be extremely handy in our analysis.

From Definition \ref{morrey-def} and H\"older's inequality we see in particular that
\be{f-avg}
 \frac{1}{|B(x,r)|} \int_{B(x,r) \cap \Omega} |f| \leq C r^{-n/p} \| f \|_{M^p_q(\Omega)};
\end{equation}
In other words, if $f \in M^p_q(\Omega)$, then $f$ has magnitude $O(r^{-n/p})$ on balls
of radius $r$, in some $L^q$-averaged sense.

We now develop some basic estimates on Morrey spaces.  All our functions here 
will be assumed to be smooth; it will not make a difference whether the functions are scalar, vector, or 2-form valued 
since we are allowing our constants $C$ to depend on $N$.  
In this section we shall also allow the constants $C$ to depend 
on the exponents $p$, $q$.

From H\"older's inequality we have
$$\| fg \|_{M^p_q(\Omega)} \leq C \|f\|_{M^{p_1}_{q_1}(\Omega)} 
\|g\|_{M^{p_2}_{q_2}(\Omega)} 
$$
for arbitrary $f,g,\Omega$,  whenever $1/p = 1/p_1 + 1/p_2$ and $1/q = 1/q_1 
+ 1/q_2$.
In particular, if $\Omega$ has finite measure, then $M^{p_1}_{q_1}(\Omega)$ 
embeds into $M^p_q(\Omega)$.  Also we have
$$
\| fg \|_{M^p_q(\Omega)} \leq C \| f\|_{M^p_q(\Omega)} \|g\|_{L^\infty(\Omega)},
$$
and that the dual of $M^p_q$ is $M^{p'}_{q'}$ when $1 < q \leq p < \infty$.  
Finally, we have the trivial observation
$$ \| \nabla^j f \|_{M^p_{q,k}(\Omega)} \leq C \| f \|_{M^p_{q,k+j}(\Omega)}.$$
We shall use the above estimates so frequently in the sequel that we shall not 
explicitly mention them again.

We now develop give analogues of standard harmonic analysis estimates for the 
Morrey space setting.

\begin{proposition}\label{smoot}
Let $1 < q \leq p < \infty$, and let $T$ be a pseudo-differential operator of 
order 0.  Then $T$ is bounded on $M^p_q(\R^n)$.

The same result holds if $T$ is replaced by the Hardy-Littlewood maximal  
operator
$$ \M u(x) := \sup_{r > 0} \frac{1}{|B(x,r)|} \int_{B(x,r)} |u|.$$
\end{proposition}

\begin{proof}
We give the proof for $T$ only, as the argument for $\M$ is completely analogous.

We need to show that
$$ r^{n(\frac{1}{p}-\frac{1}{q})} \| Tf \|_{L^q(B(x,r))} \leq C \| f 
\|_{M^p_q(\R^n)}$$
for all balls $B(x,r)$.  By scaling we may take $B(x,r) = B(0,1)$.

First suppose that $f$ is supported on $B(0,2)$.  Then the claim follows from 
the standard result that $T$ is bounded on $L^q$ (see e.g. \cite{stein:large}).  
Thus we may assume that $f$ vanishes on $B(0,2)$.  In this case we use the fact 
that the kernel $K(x,y)$ of $T$ must decay like $O(|x-y|^n)$ and a standard 
dyadic decomposition to obtain the pointwise estimate
$$ |Tf(x)| \leq C \sum_{k=0}^\infty 2^{-nk} \int_{B(0,2^k)} |f|.$$
The claim then follows from \eqref{f-avg}.
\end{proof}

As an immediate corollary of this proposition we see that any smoothing operator 
of order $k$ will map $M^p_q$ to $M^p_{q,k}$ whenever $1 < q \leq p < \infty$.

We now develop further corollaries of the above Proposition.

\begin{proposition}[Fractional integration]\label{frac}
Whenever $n/2 \leq p < n$ and $1/q = 1/p - 1/n$, we have
$$
\| u * \frac{1}{|x|^{n-1}} \|_{M^q_4(\R^n)} \leq C_{p, q} \| u 
\|_{M^{p}_2(\R^n)}.$$
\end{proposition}

\begin{proof}
We may assume that $\| u \|_{M^p_2(\R^n)} = 1$.

Let $x_0$ be an arbitrary point in $\R^n$, and use dyadic decomposition and 
H\"older to estimate
$$ |u * \frac{1}{|x|^{n-1}}(x_0)| \leq
\sum_{k \in \Z: 2^k < r} 2^{-(n-1)k} \| u \|_{L^1(B(x_0,2^k))}
+ \sum_{k \in \Z: 2^k \geq r}
2^{-(n-1)k} 2^{nk/2} \| u \|_{L^2(B(x_0,2^k))}$$
where $r > 0$ will be chosen later.

We can bound the first term by $C r \M u(x_0)$.  To control the second term, we 
use the bound $\| u \|_{L^2(B(x_0,2^k))} \leq C_p 
r^{n(\frac{1}{2}-\frac{1}{p})}$ and the assumption $p<n$ to control this by 
$r^{1-n/p}$.  Adding the two estimates together and optimizing in $r$ we thus 
have (after some algebra) the pointwise estimate
$$ |u * \frac{1}{|x|^{n-1}}(x_0)| \leq C_p \M u(x_0)^{q/p}.$$
Thus
$$ \| u * \frac{1}{|x|^{n-1}} \|_{M^q_4} \leq
C_p \| \M u \|_{M^p_{4p/q}}^{q/p}.$$
Since $p \leq n/2$ and $1/q = 1/p - 1/n$, we have $4p/q \leq 2$.  The claim then 
follows from Proposition \ref{smoot}.
\end{proof}

We now specialize our domain $\Omega$ to the unit cube $[0,1]^n$.

\begin{corollary}[Morrey-Sobolev embeddings]\label{Sobolev}
We have the estimate
\be{sobolev}
\| u \|_{M^q_4([0,1]^n)} \leq C_{p, q} \| u \|_{M^{p}_{2,1}([0,1]^n)}
\end{equation}
whenever $n/2 \leq p < n$ and $1/q = 1/p - 1/n$.
If $p$ is strictly greater than $n/2$, we also have
\be{l-infty}
\| u \|_{L^\infty([0,1]^n)} \leq C \| u \|_{M^p_{2,2}([0,1]^n)};
\end{equation}
in fact we can replace $L^\infty$ by the H\"older space $C^{0,\alpha}$ for some $0 < \alpha = \alpha(p) < 1$.
\end{corollary}

\begin{proof}
We first prove \eqref{sobolev}. From the fundamental theorem of calculus and 
polar co-ordinates around $x_0$ we have the pointwise estimate
$$ |u(x_0)| \leq C (|\nabla u| * |x|^{1-n})(x_0) + \int_{[0,1]^n} |u|.$$
The former term is acceptable by Proposition \ref{frac}.  The latter term is 
acceptable by \eqref{f-avg}.

The claim \eqref{l-infty} (and the H\"older refinement) then follows from \eqref{sobolev} and Morrey's lemma, see e.g. \cite{taylor}.  
\end{proof}

From the above Proposition and H\"older we obtain the basic estimate
\be{hybrid-y}
\| uv \|_{M^{n/2}_2([0,1]^n)} \leq C \| u \|_{M^{n/2}_{2,1}([0,1]^n)} \| 
v\|_{M^{n/2}_{2,1}([0,1]^n)}.
\end{equation}
This estimate is what allows us to control the non-linear portion $A \wedge A$ 
of the curvature by the linear portion $dA$, assuming that one is in the Coulomb 
gauge and that the curvature is small in $M^{n/2}_2$.

For technical reasons having to do with continuity arguments we shall need to also
work in slightly smoother spaces than $M^{n/2}_{2,k}$, and in particular we
shall need to investigate the continuity of the Coulomb 
gauge construction in the space $M^p_{2,1}$ for some $n/2 < p < n$.

For these spaces one can use
Corollary \ref{Sobolev}, H\"older, and the Leibnitz rule for derivatives to 
obtain the product estimates
\be{algebra}
\| u v \|_{M^p_{2,i+j-2}([0,1]^n)} \leq C_p \| u\|_{M^p_{2,i}([0,1]^n)} \| 
v\|_{M^p_{2,j}([0,1]^n)}
\hbox{ whenever } i,j,i+j-2 \in \{0,1,2\}.
\end{equation}
Also, we shall need the variant
\be{hybrid-x}
\| uv \|_{M^p_2([0,1]^n)} \leq C_p \| u \|_{M^p_{2,1}([0,1]^n)} \| 
v\|_{M^{n/2}_{2,1}([0,1]^n)}
\end{equation}
which is proven by the same methods.

We now give some standard elliptic regularity estimates for Morrey spaces.  In 
proving these estimates it is convenient to define the 
approximate fundamental solution $K$ to the Laplacian by
$$ K := c\phi/|x|^{2-n}$$
where $\phi$ be a smooth radial bump function which equals 1 on $[-2,2]^n$ and 
$c := -4 \pi^{n/2} / \Gamma(\frac{n-2}{2})$ is the constant such that $\Delta \frac{c}{|x|^{2-n}}$ is the Dirac delta.  
Observe that $\Delta K = \delta + \psi$ for some bump function $\psi$.

\begin{proposition}[Elliptic regularity]\label{elliptic-prop}
Let $1 < q < \infty$.  If $u^\alpha$, $\varphi^{\alpha\beta}$ satisfy the Hodge 
system
\begin{align}
\partial^\beta u^\alpha - \partial^\alpha u^\beta &= \varphi^{\alpha\beta}
\label{ab}\\
\partial_\alpha u^\alpha &= 0 \label{ac}
\end{align}
on $[0,1]^n$ with the boundary condition
\be{ad}
n_\alpha u^\alpha = 0
\end{equation}
on $\partial [0,1]^n$, then
\be{elliptic-2}
\| u \|_{M^q_{2,1}([0,1]^n)} \leq C_q \| \varphi \|_{M^q_2([0,1]^n)}.
\end{equation}
\end{proposition}

\begin{proof}
Extend the one-form $u^\beta$ and the two-form $\varphi^{\alpha \beta}$ to 
$\R^n$ so that they are symmetric with respect to reflections across the faces 
of $[0,1]^n$.  Observe from \eqref{ab}, \eqref{ac}, \eqref{ad} that \eqref{ab}, 
\eqref{ac} in fact hold on all of $\R^n$ in the sense of distributions.

Contracting \eqref{ab} against $\partial_\beta$ and using \eqref{ac}, we obtain
$$ \Delta u^\alpha = \partial_\beta \varphi^{\alpha\beta}.$$
Motivated by this, we define
$$ \tilde u^\alpha := \partial_\beta \varphi^{\alpha\beta} * K.$$
We thus see that $\Delta (u^\alpha - \tilde u^\alpha)$ is smooth on $[0,1]^n$, 
and that $n_\alpha (u^\alpha - \tilde u^\alpha) = 0$ on $\partial [0,1]^n$.  
From this and standard elliptic regularity one sees that the contribution of $u - \tilde u$ is acceptable.  To deal with the contribution of $\tilde u$ we just 
observe that convolution with $\partial_\beta K$ is a standard smoothing 
operator of order 1, and use Proposition \ref{smoot}.
\end{proof}

\begin{proposition}[Neumann problem regularity]\label{neumann-prop}
Let $1 < q < \infty$.  Let $u, f, g_\beta \in C^\infty([0,1]^n)$ be such that
\be{laplace}
\Delta u = f
\end{equation}
on $(0,1)^n$ and
\be{boundary}
n^\beta \partial_\beta u = n^\beta g_\beta
\end{equation}
on $\partial [0,1]^n$, and we have the normalization
\be{mean}
\int_{[0,1]^n} u = 0.
\end{equation}
Then we have
\be{neumann}
\| u \|_{M^q_{2,2}([0,1]^n)} \leq C \| f \|_{M^q_2([0,1]^n)} + C \| g 
\|_{M^q_{2,1}([0,1]^n)}.
\end{equation}
\end{proposition}

\begin{proof}
We first prove the claim when $f=0$.  From the Sobolev trace lemma we observe 
that $g_\beta$ is in $L^2$ on hyperplanes.

For each $\beta$, we extend $g_\beta$ to the domain $\{ x: x_\beta \in [0,1] \}$
by requiring $g_\beta$ to be symmetric with respect to reflections along the 
faces of $[0,1]^n$ perpendicular to $e_\beta$.
Define $\tilde u$ on $[0,1]^n$ by
$$ \tilde u(x) := \sum_{\beta = 1}^n \int_{y_\beta = 1} K(x-y) g_\beta(y)\ dy - 
\int_{y_\beta = 0} K(x-y) g_\beta(y)\ dy$$

Since $\Delta K = \delta + \psi$ we have
$$ \Delta \tilde u(x) = \sum_{\beta = 1}^n \int_{y_\beta = 1} \psi(x-y) 
g_\beta(y)\ dy - \int_{y_\beta = 0} \psi(x-y) g_\beta(y)\ dy$$
on $(0,1)^n$.  In particular, $\Delta \tilde u$ is smooth, with a norm 
controlled by the $M^q_{2,1}([0,1]^n)$ norm of $g$.  From the Plemelj formulae 
and symmetry we also see that $n^\beta \partial_\beta \tilde u = n^\beta 
g_\beta$ on $\partial [0,1]^n$.  Thus, it remains only to show that $\tilde u$ 
is in $M^q_{2,2}([0,1]^n)$, since the difference $u - \tilde u$ can be 
controlled easily by the classical Neumann theory.

By symmetry it suffices to show that
$$ \|\int_{y_n = 0} K(x-y) g(y)\ dy \|_{M^q_{2,2}([0,1]^n)} \leq C \| g 
\|_{M^q_{2,1}(\R^n)}$$
for any function $g$ on $\R^n$.  By applying a cutoff we may assume that $g(y)$ 
is supported on the region $|y_n| \leq 1$.  We may then invoke the Fundamental 
theorem of calculus and polar co-ordinates to write
$$ g = \nabla g * L$$
for some (covector-valued) kernel $L$ supported on the cone $\{ x: |x| \leq x_n 
\leq C \}$ which behaves like $|x|^{1-n}$.  It thus suffices to show that
$$ \|\nabla_x^k \int \int_{y_n = 0} K(x-y) L(y-z) G(z)\ dy\ dz 
\|_{M^q_{2,2}([0,1]^n)} \leq C \| G \|_{M^q_2(\R^n)}$$
for all $G$ and $k=0,1,2$.

The expression inside the norm vanishes for $|y_n| \geq C$.  From this and the 
fundamental theorem of calculus we see that we need only prove the above 
estimate for $k=2$.

A computation shows
$$ |\nabla_x^2 \int_{y_n = 0} K(x-y) L(y-z)\ dy| \leq x_n^{-n} (1 + 
|z-x|/x_n)^{-n-1}$$
for all $x$, $z$, which implies the pointwise estimate
$$ |\nabla_x^k \int \int_{y_n = 0} K(x-y) L(y-z) G(z)\ dy\ dz|
\leq C \M G(x).$$
The claim then follows from Proposition \ref{smoot}.

We now consider the general case when $f$ is not necessarily 0.   Define
$$ v := u - (f\chi_{[0,1]^n}) * K.$$
Then
$$ \partial_\alpha \partial^\alpha v = 0$$
on $[0,1]^n$ and
$$ n^\beta \partial_\beta v = n^\beta (g_\beta - f\chi_{[0,1]^n} * 
\partial_\beta K).$$
The operation of convolution with $K$ is smoothing of order 2, thus
$$ \| (f\chi_{[0,1]^n}) * K \|_{M^q_{2,2}} \leq C \| f \|_{M^q_2} $$
and \eqref{neumann} then follows from the $f=0$ estimate applied to $v$.
\end{proof}

\begin{lemma}[Interior regularity]\label{ellip-lemma}
Let $B(x,r)$ be a ball, and let $0 < \theta \ll 1$.  Then we have
$$ \| u \|_{M^{n/2}_{2,1}(B(x,\theta r))}
\leq C \| \Delta u \|_{M^{3n/2}_{4/3}(B(x,r))}
+ C \theta^2 \| u \|_{M^{n/2}_{2,1}(B(x,r))}$$
whenever $u$ is such that the right-hand side makes sense.
\end{lemma}

\begin{proof}
By the usual limiting argument (using the ellipticity of $\Delta$) we may assume 
that $u$ is smooth.  We may rescale $B(x,r)$ to be $B(0,1)$.

Let $\eta$ be a bump function adapted to $B(0,1)$ which equals 1 on $B(0,1/2)$.  
We split
$$ u = K * (\eta \Delta u) + (u - K * (\eta \Delta u)).$$
Consider the latter term.  This is harmonic on $B(0,1/2)$, so by interior 
regularity we have
$$ \| \nabla^j (u - K * (\eta \Delta u)) \|_{L^\infty(B(0,\theta))}
\leq C \| u - K * (\eta \Delta u) \|_{L^1(B(0,1/4))}$$
for $j=0,1$.  The right-hand side can be easily bounded by
$$ C \| \Delta u \|_{M^{3n/2}_{4/3}(B(0,1))} + C \| u 
\|_{M^{n/2}_{2,1}(B(0,1))}.$$
Since for any $F$ we have
$$ \| F \|_{M^{n/2}_2(B(0,\theta))}
\leq C \theta^2
\| F \|_{L^\infty(B(0,\theta))}
$$
we thus see that the contribution of this term is acceptable.

To control the first term it suffices to show the global estimate
$$ \| \nabla^j K * f \|_{M^{n/2}_2(\R^n)} \leq C \| f \|_{M^{3n/2}_{4/3}(\R^n)}$$
for test functions $f$ and $j=0,1$.  But this follows from the dual of Proposition 
\ref{frac} (with $p=n/2$ and $q=n$), since $\nabla^j K$ is pointwise bounded by 
$C/|x|^{n-1}$.
\end{proof}

\section{Overview of proof of Theorem \ref{abstract3}}

We now give the proof of Theorem \ref{abstract3}, modulo some Propositions which we will prove in later sections.

Fix $A$ to be an admissible Yang-Mills connection obeying the assumptions in Theorem \ref{abstract3}.
From the analysis in \cite{tian} we have the following properties:

\begin{proposition}\label{A-ass}
Let $A$ be a stationary admissible Yang-Mills connection and $S$ be its
singular set. Then the curvature $F(A)$ 
obeys the Morrey norm estimate
\be{curv-small}
\| F(A) \|_{M^{n/2}_2([0,1]^n)} \leq C\eps
\end{equation}
and the pointwise bounds
\be{f-point}
|F(A)(x)| \le
\frac{C\eps}{\rho(x)^2} \hbox{ for all } x\in [0,1]^n\backslash S,
\end{equation}
where $\rho(x)$ is the distance function
$$\rho(x) := \dist(x,S).$$
\end{proposition}

(\ref{curv-small}) follows from the monotonicity for stationary
Yang-Mills connections. (\ref{f-point}) is obtained by applying
curvature estimates to the smooth connection $A$ outside $S$
and using (\ref{curv-small}).
The singular set is extremely small, having codimension at least 4.  This high codimension shall be
crucial in our arguments, as it allows various path and surface integrals to be generically well-defined.  Note that we have no control over the topology or regularity of $S$ (other than that $S$ is compact), however we will use averaging arguments to get around this difficulty.

We remark that the next few steps of the argument (Lemma \ref{size}, Proposition \ref{inductive-gauge}, Proposition \ref{m-approx}, Theorem \ref{abstract})
do not use the Yang-Mills equation \eqref{ym} directly; instead, they proceed from the conclusions in Proposition \ref{A-ass}, which of course hold for a more general class of connections than the stationary admissible Yang-Mills connections. The Yang-Mills equation only re-appears much later in the argument, in Lemma \ref{delta-a}.

The function $F(A)$ is defined a.e. on $[0,1]^n$; we extend it by zero outside of this unit cube.  
From \eqref{curv-small} and \eqref{f-avg} we observe in particular that
\be{F-avg}
\int_{B(x,r)} |F(A)| \leq C \eps r^{n-2}
\end{equation}
for all balls $B(x,r)$.

To visualize a function $F(A)$ which obeys \eqref{curv-small}, \eqref{f-point}, one model to
keep in mind is when $S$ is a smooth $n-4$-dimensional surface, and $|F(A)(x)| \sim \eps / \rho(x)^2$.
(This function is not quite in $L^2$, and so \eqref{curv-small} is not quite obeyed, but this can
be rectified by adding a logarithmic decay factor).  Using the heuristic that the connection $A$
is like an antiderivative of the curvature $F(A)$, one then expects to be able to place $A$
in a gauge $\sigma(A)$ so that $|\sigma(A)(x)| \leq C \eps / \rho(x)$.  Another formulation of this 
heuristic is that one expects to be able to find a gauge $\sigma(A)$ obeying the pointwise bounds
\be{sfa}
|\sigma(A)(x)| \leq C \int_{[0,1]^n} \frac{|F(A)(y)|}{|x-y|^{n-1}}\ dy,
\end{equation}
since the right-hand side is essentially the fractional integral $|\nabla|^{-1} |F(A)|$ of $|F(A)|$.  
Such a gauge is easy to obtain in the case when $G$ is abelian, since one can use Hodge theory
to find a gauge $\sigma(A)$ which is essentially equal to $\Delta^{-1} d_* F \approx \nabla^{-1} F$.  An essentially equivalent approach in the abelian case is to take various radial gauges (where $\sigma(A)(x) \cdot (x-x_0) = 0$ for all $x \in [0,1]^n$) and average over all choices of origin $x_0 \in [0,1]^n$ to obtain a connection obeying the bounds \eqref{sfa}.  Note that this averaging will eliminate the difficulty that the rays in the radial gauge occasionally pass through the singular set $S$, since this set has such high codimension.

We do not know how to achieve bounds of the form \eqref{sfa} in the non-abelian case, however we can develop a reasonably good rigorous substitute for this heuristic, which we now discuss.  The idea is to inductively construct a successive sequence of partial gauges which obey reasonable connection bounds and which advance closer and closer to the singular set $S$ as the induction progresses.

We shall need a fixed exponent $0 < \kappa < 1$; for sake of concreteness we set $\kappa := 1/2$. We 
define\footnote{We apologize for the artificiality of the quantity $Q(x)$ and the companion quantity $T_m(x)$ 
defined in \eqref{T_m-def}; these choices were obtained after much trial and error.  One needs $Q$ small enough 
that one has the density bounds in Lemma \ref{size} and the clustering bounds in Lemma \ref{f-cluster-m}, but $Q$ 
large enough that the errors arising from the truncation argument in Proposition \ref{m-approx} are manageable.  For purposes of dimensional analysis, $Q$ has units of $length^{\kappa-1}$, $T_m$ has units of $length^{-1}$, while $\rho(x)$ and radii such as $r$ or $R_m = C D^{-m}$ have the units of $length$.  This dimensional analysis can serve to explain many of the strange powers of $r$ or $R_m$ which appear in many of the estimates in the next few sections.} the quantity
 $Q(x)$ on $\R^n$ by
\be{Q-def}
Q(x) := \sup_{0 < r \leq \diam([0,1]^n)} r^{-n/2+1+\kappa} (\int_{B(x,r)} |F(A)(y)|^2\ dy)^{1/2};
\end{equation}
this expression is like the Hardy-Littlewood maximal  function $\M F(A)$ but with an additional decay factor of
 $r^{1 + \kappa}$.
In the model case where $S$ is a smooth $n-4$-dimensional surface and $F$ is comparable to $\eps / \rho(x)^2$, the
 quantity $Q(x)$ is comparable to $\eps / \rho(x)^{1-\kappa}$; the reader may find this model case helpful to keep 
in mind in what follows.

Let $D \gg 1$ be a large number depending only on $n$, $G$, $N$ to be chosen later; we will assume that $\eps$ is sufficiently small depending on $D$.  This quantity $D$ shall be our dyadic base, i.e. we will use powers $D^m$ of $D$ to define dyadic scales, as opposed to the more usual powers of two.

We define the domains $\Omega_m$ for $m = 1, 2, \ldots$ by
\be{omega-def}
\Omega_m := \{ x \in [0,1]^n \backslash S: Q(x) < \eps D^{(1-\kappa)m} \}.
\end{equation}
Thus in the model case, $\Omega_m$ is essentially the region where $\rho(x) \geq D^{-m}$.

In general, the $\Omega_m$ are a non-decreasing collection of open sets.  As $m$ gets larger, the set $\Omega_m$ 
fills out an increasingly large portion of the unit cube, as the following lemma indicates:

\begin{lemma}[$\Omega_m$ is dense at scale $R_m$]\label{size} For any $m \geq 1$, define the radius
\be{R_m-def}
R_m := C D^{-m}
\end{equation}
where $C$ is a sufficiently large constant. Then we have
\be{omega-contains}
\{ x \in [0,1]^n \backslash S: \rho(x) \geq R_m \}
 \subseteq \Omega_m\end{equation}
Furthermore, we have
\be{large-intersection}
| B(x,r) \cap \Omega_m | \geq C^{-1} r^n \hbox{ whenever } r \geq R_m, x \in [0,1]^n.
\end{equation}
\end{lemma}

\begin{proof}
We first prove \eqref{omega-contains}.  Let $x \in [0,1]^n$ be such that $\rho(x) \geq R_m$.  We have to show that
$$ r^{-n/2+1+\kappa} (\int_{B(x,r)} |F(A)(y)|^2\ dy)^{1/2} \leq \eps D^{(1-\kappa) m}$$
for all $0 < r \leq \diam([0,1])$.

First suppose that $r < \rho(x)/2$.  Then by \eqref{f-point} we can bound the left-hand side by
$$ r^{-n/2+1+\kappa} r^{n/2} \eps / \rho^2(x) \leq C \eps \rho(x)^{-1+\kappa} \leq C \eps R_m^{-1+\kappa},$$
which will be acceptable if the constant in \eqref{R_m-def} is large enough.

Now suppose that $r > \rho(x)/2$.  Then by \eqref{curv-small} we can bound the left-hand side by
$$ r^{-n/2+1+\kappa} \eps r^{n/2-2} \leq C \eps \rho(x)^{-1+\kappa} \leq C \eps R_m^{-1+\kappa}$$
which is again acceptable.

Now we prove \eqref{large-intersection}.  Fix $B(x,r)$.  For each $x'$ in the compact set
$\overline{B(x,r)} \backslash \Omega_m$
we see, from \eqref{omega-def}, \eqref{Q-def}, that there exists a radius $0 < r(x') \leq C$ such that
$$ r(x')^{-n/2+1+\kappa} (\int_{B(x',r(x'))} |F(A)(y)|^2\ dy)^{1/2} \geq \eps D^{(1-\kappa) m}.$$
From \eqref{curv-small} we must have $r(x') \leq C D^{-m}$.  In particular we have $B(x',r(x')) \subseteq
B(x, 2r)$ if the constant in \eqref{R_m-def} is large enough.

The balls $B(x',r(x'))$ clearly cover the compact set $\overline{B(x,r)} \cap ([0,1]^n \backslash \Omega_m)$,
so in particular there is a finite sub-cover of this set by these balls.  By the Vitali covering lemma there
thus exists a finite sub-collection $B(x_j, r(x_j))$ of balls which are disjoint and such that
$B(x_j, 5r(x_j))$ covers $\overline{B(x,r)} \cap ([0,1]^n \backslash \Omega_m)$.  In particular we have
$$ |B(x,r) \cap [0,1]^n \backslash \Omega_m| \leq C \sum_j r(x_j)^n.$$
On the other hand, by construction we have
$$ \int_{B(x_j, r(x_j))} |F(A)(y)|^2\ dy
\geq \eps^2 D^{2(1-\kappa) m} r(x')^{n-2-2\kappa}
\geq C^{-1} \eps^2 D^{4m} r(x')^n$$
since $r(x') \leq C D^{-m}$.  Thus we have
$$ |B(x,r) \cap [0,1]^n \backslash \Omega_m| \leq \sum_j \frac{C}{\eps^2 D^{4m}} \int_{B(x_j,r(x_j))} |F(A)(y)|^2\ dy;$$
since the balls $B(x_j, r(x_j))$ are disjoint and lie in $B(x,2r)$ we thus have
$$ |B(x,r) \cap [0,1]^n \backslash \Omega_m| \leq \frac{C}{\eps^2 D^{4m}} \int_{B(x,2r)} |F(A)(y)|^2\ dy;$$
from \eqref{curv-small} we thus have
\bas |B(x,r) \cap [0,1]^n \backslash \Omega_m| \leq C D^{4m} r^{n-4} &\leq C |B(x,r) \cap [0,1]^n| (r D^m)^{-4} \\
&\leq C |B(x,r) \cap [0,1]^n| (r/R_m)^{-4}.
\end{align*}
If the constant in \eqref{R_m-def} is sufficiently large, the claim \eqref{large-intersection} follows.
\end{proof}

We are now ready the precise analogue of the heuristic \eqref{sfa}.

\begin{proposition}\label{inductive-gauge}  Let $A$, $S$ be as in Proposition \ref{A-ass}.
Then, if $D$ is sufficiently large and $\eps$ is sufficiently small depending on $D$, for every $m \geq 1$ 
there exists a gauge transform $\sigma_m$ on $\Omega_m$ which obeys the bounds
\be{sigma_m-bound-alt}
|\sigma_m(A)(x)| \leq C(D) T_m(x).
\end{equation}
for all $x \in \Omega_m$, where $T_m(x)$ is the quantity
\be{T_m-def}
T_m(x) := \int_{[0,1]^n} (1 + \frac{|y-x|}{R_m})^{\kappa/2} 
\frac{|F(A)(y)|}{|x-y|^{n-1}}\ dy.
\end{equation}
\end{proposition}

Using \eqref{F-avg} one can crudely bound $T_m(x)$ by $O(\eps/R_m)$, but we will need the more precise structure of 
\eqref{T_m-def} in our truncation analysis later in this section.  However, this crude bound already
shows that $\sigma_m(A)$ is bounded on $\Omega_m$, and hence (by \eqref{e-def} and the smoothness of $A$) that 
$\sigma_m$ is locally Lipschitz on $\Omega_m$, although our Lipschitz bound of course depends on $m$.

We do not assert here that $\sigma_m$ is smooth; however when we use this Proposition later in this section we will be 
able to regularize $\sigma_m$ with little difficulty.

The construction of the gauges $\sigma_m$ will be inductive, with the gauge $\sigma_{m+1}$ obtained from $\sigma_m$ by
integrating the connection $A$ along curves, using the curvature 
bound \eqref{curv-small} to compare integrals along two different curves, and averaging over a family of curves; we 
give the construction and prove the Proposition in Sections \ref{integrating-sec}-\ref{inductive-sec}.  It may be 
possible to obtain a result like this
more directly, perhaps by using the finer structure of the Yang-Mills equation, but we were unable to do so.

Using Proposition \ref{inductive-gauge} and a truncation argument, we will prove the following approximation theorem:

\begin{proposition}[Approximation by smooth connections]\label{m-approx}  
For any $m > 0$, there exists a smooth connection $A_m$ on $[0,1]^n$ which is 
gauge equivalent to $A$ on the open set
\be{open}
\{ x \in [0,1]^n: \rho(x) \geq 20 R_m \}
\end{equation}
and which obeys the smallness condition
$$ \| F(A_m) \|_{M^{n/2}_2([0,1]^n)} \leq C(D) \eps.$$
\end{proposition}

We shall prove this Proposition in Section \ref{approx-sec}.

The next step is to place the smooth connections $A_m$ in a Coulomb gauge.

\begin{definition}\label{coulomb-def}
A connection $A$ is said to be a \emph{Coulomb gauge on $\Omega$} if it 
satisfies the condition
$d_* A = 0$
on the interior of $\Omega$, and
$A \cdot n = 0$
on the boundary $\partial \Omega$.
\end{definition}

From elliptic theory we expect Coulomb gauges to be quite
regular; specifically, we expect $A$ to have one more derivative
of regularity than $F(A)$. The question then arises: given an
arbitrary connection $A$, under what conditions can we find a
gauge equivalent Coulomb gauge $A_{coulomb}$ which has a one more
derivative of regularity than $F(A)$?

In \cite{uhlenbeck} this problem was considered assuming that the curvature was 
in $L^{n/2}$.  In our applications we need to replace this space by the slightly 
larger space $M^{n/2}_2$.  Using the heuristic that a connection requires one more derivative than the 
curvature, and a gauge transform requires two more derivatives, we thus hope to 
place connections and gauge transforms in $M^{n/2}_{2,1}$ and $M^{n/2}_{2,2}$ 
respectively.

For any $\eps > 0$, let $\U_{\eps}(\Omega)$ denote the set of all smooth 
connections on $\Omega$ which satisfy the bound
\be{f-small}
\| F(A) \|_{M^{n/2}_2(\Omega)} \leq \eps.
\end{equation}
From \eqref{gauge-eq} we observe that this space is invariant under gauge 
transformations.

In Section \ref{setup-sec} we prove the following generalization of Uhlenbeck's lemma \cite{uhlenbeck}:

\begin{theorem}[Small curvature allows a Coulomb gauge]\label{abstract} If $0 < \eps \ll 1$ is sufficiently small, then 
every connection $A$ in $\U_\eps([0,1]^n)$ is gauge equivalent (via a smooth gauge transformation $\sigma$) to a smooth 
Coulomb 
gauge $A_{coulomb}$ which obeys the bound
\be{bound1}
\| A_{coulomb} \|_{M^{n/2}_{2,1}([0,1]^n)} \leq C \| F(A) 
\|_{M^{n/2}_2([0,1]^n)}.
\end{equation}
\end{theorem}

This theorem has also been independently obtained by Riviere and Meyer \cite{riviere}.  Note that Theorem \ref{abstract}
does not require $A$ to be Yang-Mills, but it does require $A$ to be smooth.

We apply Theorem \ref{abstract} to the smooth gauges $A_m$ constructed in Proposition \ref{m-approx}.
We may thus find (if $\eps$ is sufficiently small depending on $C_0$, $D$) a Coulomb gauge
$A_{m,coulomb} := \sigma^{coulomb}_m(A_m)$ such that
\be{j-loc} \| A_{m,coulomb} \|_{M^{n/2}_{2,1}([0,1]^n)} \leq C \eps.
\end{equation}
By weak compactness, we may thus find a sequence $m_j \to 0$ such that
$A_{m_j,coulomb}$ converges weakly in $M^{n/2}_{2,1}([0,1]^n)$ to a function 
$A_{coulomb}$ such that
\be{c-bound}
\| A_{coulomb} \|_{M^{n/2}_{2,1}([0,1]^n)} \leq C \eps.
\end{equation}
By taking weak limits of the Coulomb gauges $A_{m_j,coulomb}$ we thus see that 
$A_{coulomb}$ is also a Coulomb gauge\footnote{Observe from \eqref{c-bound} that the $A_{m_j,coulomb}$ are uniformly 
in the Sobolev space $W^{1,2}$, and so by Rellich embedding and the Sobolev trace lemma they converge strongly in $L^2$
on the boundary of $[0,1]^n$, thus we may legitimately take limits of the boundary condition
$A_{m_j,coulomb} \cdot n = 0$.}.

Let $K$ be an arbitrary compact subset of $[0,1]^n \backslash S$.  Then for $j$ sufficiently large, we have 
that $R_{m_j} \ll \dist(K,S)$, so that $A_{m_j,coulomb}$ is gauge equivalent to $A$ on 
$K$.  Thus we may find a gauge $\overline \sigma_j$ on $K$ such that
\be{limit}
d \overline \sigma_j = \overline \sigma_j A - A_{m_j, coulomb} \overline \sigma_j.
\end{equation}
Since $A$, $A_{m_j,coulomb}$ are in $L^2$ uniformly in $j$, we thus see that the
$\nabla \overline \sigma_j$ are in $L^2$ uniformly in $j$.
By passing to a subsequence if necessary and using Rellich embedding we see that 
$\overline \sigma_j$ converges strongly in $L^q$ for some $q > 2$ to some limiting function $\sigma \in L^q$.  
In particular we see 
that $\sigma$ takes values in $G$.  By taking limits of \eqref{limit} we thus see that
$$ d\sigma = \sigma A - A_{coulomb} \sigma$$
and in particular that $\sigma \in W^{r,1}$ for some $r>1$.  Thus $A$ and $A_{coulomb}$ are gauge equivalent on $K$, 
and hence (since $K$ was arbitrary)
on $[0,1]^n \backslash S$.  

To summarize, we have found a gauge $A_{coulomb} = \sigma(A)$ on $[0,1]^n \backslash S$ which is a Coulomb gauge and 
obeys the smallness assumption \eqref{c-bound}.  To complete the proof of Theorem \ref{abstract3}, it will suffice
 to show that $A_{coulomb}$ extends to be smooth on all of $[0,1]^n$.

We shall use standard elliptic regularity techniques.  The first step is to obtain some bounds on  $\Delta A_{coulomb}$,
 which presently is only defined in the sense of distributions.

\begin{lemma}\label{delta-a}  The distribution $\Delta A_{coulomb}$ is a locally integrable function; in fact it lies
in the space $M^{3n/2}_{4/3}([0,1]^n)$ and obeys the pointwise estimate
\be{jerry}
|\Delta A_{coulomb}| \leq C | A_{coulomb} | |\nabla A_{coulomb}|
+ C |A_{coulomb}|^3
\end{equation}
almost everywhere in $[0,1]^n$.
\end{lemma}

\begin{proof}
First we work in $[0,1]^n \backslash S$.  On this set $A_{coulomb}$ is a Yang-Mills connection in the Coulomb gauge, 
so by \eqref{ym} and the Coulomb gauge condition we have
\bas
\Delta A_{coulomb} &= d_* dA_{coulomb} + d d_* A_{coulomb} \\
&= d_* d A_{coulomb} \\
&= d_*(F(A_{coulomb})) - d_* (A_{coulomb} \wedge A_{coulomb}) \\
&=  \ast [A_{coulomb}, \ast F(A_{coulomb})] - d_* (A_{coulomb} \wedge A_{coulomb}) \\
&= O(|F(A_{coulomb})| |A_{coulomb}| + |\nabla A_{coulomb}| |A_{coulomb}|) \\
&= O(|\nabla A_{coulomb}| |A_{coulomb}| + |A_{coulomb}|^3)
\end{align*}
on $[0,1]^n \backslash S$.  Note that $F(A_{coulomb})$ is in $M^{n/2}_2$ and $A_{coulomb}$ is in $M^n_4$ 
(by \eqref{c-bound}, \eqref{sobolev}), so the above computations are justified in the sense of distributions.  
In particular we see that \eqref{jerry} holds a.e. in $[0,1]^n \backslash S$, which by \eqref{sobolev}, 
\eqref{c-bound}, and H\"older implies that $\Delta A_{coulomb}$, \emph{when restricted to $[0,1]^n \backslash S$}, 
is locally integrable and lies in the Morrey space $M^{3n/2}_{4/3}([0,1]^n \backslash S)$.

We are almost done; however, we still have to exclude the technical possibility that the distribution 
$\Delta A_{coulomb}$ has a singular component on the set $S$. Fortunately, the high codimension of $S$ will prevent
this from happening, as $S$ is too small to support singularities with the required regularity.  

To avoid confusion, let us use 
$\Delta A_{coulomb}|_{[0,1]^n \backslash S}$ to denote the (classical) Laplacian of $A_{coulomb}$ outside of $S$; 
we have already shown that $\Delta A_{coulomb}|_{[0,1]^n \backslash S}$ obeys \eqref{jerry} and the Morrey 
space bounds.  It remains to show that $\Delta A_{coulomb} = \Delta A_{coulomb}|_{[0,1]^n \backslash S}$ in 
the sense of distributions.  Fortunately \eqref{c-bound} will provide enough regularity on $A_{coulomb}$ to achieve this.

We turn to the details.  Let $0 < \delta \ll 1$ be a small parameter (which we will eventually send to zero).  Since
$S$ has Hausdorff dimension at most $n-4$, we can find a finite number of balls $B(x_j, r_j)$ with $0 < r_j < \delta$ 
which cover $S$ and obey the bound
\be{rj-bound}
\sum_j r_j^{n-3} \leq C \delta^{1/2};
\end{equation}
indeed one could replace the exponent $1/2$ by any exponent between 0 and 1.

For each ball $B(x_j,r_j)$, let $\eta_j$ be a non-negative bump function adapted to $B(x_j,2r_j)$ which equals 1 on 
$B(x_j,r_j)$, and define
$$ \eta := \sup_j \eta_j.$$
Thus $\eta$ is a Lipschitz function supported on $\bigcup_j B(x_j,2r_j)$ which equals 1 on $\bigcup_j B(x_j,r_j)$.  
We have the easily verified pointwise bound
\be{eta-bound}
|\nabla \eta(x)| \leq C \sum_j r_j^{-1} \chi_{B(x_j,2r_j)}.
\end{equation}

From \eqref{c-bound} we see that the function $(1-\eta)A_{coulomb}$ converges weakly to $A_{coulomb}$ as $\delta \to 0$.  
In particular, $\Delta( (1-\eta)A_{coulomb} )$ converges weakly to $\Delta A_{coulomb}$ in the sense of distributions.  
On the other hand, $(1-\eta) \Delta A_{coulomb}$ is supported on $[0,1]^n \backslash S$ and so converges weakly to the 
$M^{3n/2}_{4/3}$ function $\Delta A_{coulomb}|_{[0,1]^n \backslash S}$.  To conclude the argument we have to show the 
commutator estimate 
\be{delt}
\Delta( (1-\eta) A_{coulomb} ) - (1-\eta) \Delta A_{coulomb} \rightharpoonup 0
\end{equation}
as $\delta \to 0$.

The left-hand side of \eqref{delt} is equal to 
$$ -\nabla \eta \cdot \nabla A_{coulomb} - \nabla \cdot (\nabla \eta A_{coulomb})$$
so it will suffice to show the strong $L^1$ convergence
$$ \| \nabla \eta \cdot \nabla A_{coulomb}\|_{L^1} + \| \nabla \eta A_{coulomb} \|_{L^1} \to 0 \hbox{ as } \delta \to 0.$$
By \eqref{eta-bound} the left-hand side is bounded by
$$ C \sum_j r_j^{-1} \int_{B(x_j, 2r_j)} |\nabla A_{coulomb}| + |A_{coulomb}|.$$
By \eqref{c-bound} and H\"older, we can bound this by
$$ C \sum_j r_j^{-1} \eps r_j^{n-2},$$
and the claim follows by \eqref{rj-bound}.
\end{proof}

We now use Lemma \ref{delta-a} to show the interior decay estimate 
\be{interior} \| A_{coulomb}
\|_{M^{n/2}_{2,1}(B(x,\theta r))} \leq a \| A_{coulomb}
\|_{M^{n/2}_{2,1}(B(x,r))}
\end{equation}
for all balls $B(x,r) \subset (0,1)^n$, where $0 < \theta \ll 1$ and $0 < a < 1$ 
are absolute constants to be chosen later.

By interior regularity (Lemma \ref{ellip-lemma}) we have
$$
\| A_{coulomb} \|_{M^{n/2}_{2,1}(B(x,\theta r))}
\leq C \| \Delta A_{coulomb} \|_{M^{3n/2}_{4/3}(B(x,r))}
+ C \theta^2 \| A_{coulomb} \|_{M^{n/2}_{2,1}(B(x,r))}.
$$
Applying \eqref{jerry} and using \eqref{sobolev} we see that
\bas
\| A_{coulomb} \|_{M^{n/2}_{2,1}(B(x,\theta r))} \leq &
C (\| A_{coulomb} \|_{M^{n/2}_{2,1}(B(x, r))}^2
+ \| A_{coulomb} \|_{M^{n/2}_{2,1}(B(x, r))}^3)\\
&+ C \theta^2 \| A_{coulomb} \|_{M^{n/2}_{2,1}(B(x,r))}.
\end{align*}
By \eqref{c-bound} we thus have
$$
\| A_{coulomb} \|_{M^{n/2}_{2,1}(B(x,\theta r))} \leq
(C \theta^2 + C\eps) \| A_{coulomb} \|_{M^{n/2}_{2,1}(B(x,r))}$$
which gives the desired estimate \eqref{interior} if $\theta$ and $\eps$ are sufficiently small.

Iterating \eqref{interior} we see that $A_{coulomb}$ is locally
in $M^{n/2+\delta}_{2,1}$ for some absolute constant $\delta >
0$.  Applying \eqref{sobolev}, \eqref{jerry}, and elliptic
estimates (such as variants of Lemma \ref{ellip-lemma}) we may
bootstrap the value of $\delta$ in the usual fashion to obtain
smoothness of $A_{coulomb}$ on $(0,1)^n$. We omit the details as
they are rather standard.  

The proof of Theorem \ref{abstract3} is thus complete as soon as we prove
Proposition \ref{inductive-gauge}, Proposition \ref{m-approx}, and 
Theorem \ref{abstract}.  This will be achieved in Section \ref{integrating-sec}-\ref{inductive-sec}, Section \ref{approx-sec}, and Section \ref{setup-sec}.
We remark that the proofs of these results are quite distinct, and can be read independently 
of each other.

\section{Integrating connections, curvature, and averaging arguments}\label{integrating-sec}

We now begin the proof of Proposition \ref{inductive-gauge}.  In this section we shall focus on
developing the machinery needed to prove this Proposition; more precisely, we set up some notation
for integrating connections along paths and loops, and comparing the latter with
integrals of curvature.  We also need some machinery for averaging functions on the Lie group $G$ to obtain another
function on $G$.

Throughout this section we assume that $A$ obeys the properties in Proposition \ref{A-ass}, and in particular is smooth
away from the set $S$.

If $x_0, x_1, x_2 \in [0,1]^n$, we use $\Delta(x_0,x_1,x_2)$ to denote the solid triangle with these three vertices
(or equivalently, the convex hull of $\{x_0,x_1,x_2\}$);
this is a two-dimensional surface with boundary, and so we can integrate on it using two-dimensional
Hausdorff measure $d\H^2$.

If $x_0, x_1 \in [0,1]^n$, we define $[x_0 \to x_1]$ to be the path $t \mapsto x_0 + t(x_1 - x_0)$, i.e. the oriented line segment 
from $x_0$ to $x_1$.  We use $[x_0 \to x_1 \to x_2]$ as short-hand for the concatenated path 
$[x_0 \to x_1] + [x_1 \to x_2]$; thus for instance the triangular loop
$[x_0 \to x_1 \to x_2 \to x_0]$ traverses the boundary of
$\Delta(x_0, x_1, x_2)$.

Let $[x_0 \to x_1]$ be a line segment which does not intersect the singular set $S$.  We define the group element 
$A[[x_0 \to x_1]] \in G$ by the PDE
\be{A-def}
\begin{split}
A[[x_0 \to x_0]] &= 1_G\\
(x_1 - x_0) \cdot \nabla_{x_1} A[[x_0 \to x_1]] &= A[[x_0 \to x_1]] ((x_1 - x_0) \cdot A(x_1)),
\end{split}
\end{equation}
where $1_G$ is the identity element of $G$.

Since $A$ is smooth on $[x_0 \to x_1] \subset [0,1]^n \backslash S$, we see from the Picard existence theorem that $A[[x_0 \to x_1]]$ is well-defined and takes values in $G$, indeed this quantity is essentially the radial gauge from
$x_0$ evaluated at $x_1$.  In the special case where $G$ is abelian, we have the explicit formula
$$ A[[x_0 \to x_1]] = \exp(\int_{[x_0,x_1]} A(y)\ d\H^1(y)) = \exp(\int_0^1 (x_1-x_0) \cdot A(x_0 + t(x_1-x_0))\ dt)$$
but in general no such explicit formula is available. Nevertheless, it is helpful to think of $A[[x_0 \to x_1]]$
as some sort of non-abelian integral of $A$ from $x_0$ to $x_1$.

We define
$$ A[[x_0 \to x_1 \to x_2]] := A[[x_0 \to x_1]] A[[x_1 \to x_2]],$$
whenever $[x_0 \to x_1 \to x_2]$ does not intersect $S$, and similarly for more complicated polygonal paths.
The expression $A[[x_0 \to \ldots \to x_n]]$ can be regarded as the transport of the identity group element $1_G$
along the path $[x_0 \to \ldots \to x_n]$ by the connection $A$; in particular in the case of a loop $x_n = x_0$,
this element represents the monodromy of the connection along the loop.  In particular for a small triangular loop
we have
\be{f-tri}
 A[[x_0 \to x_0 + \eps v_1 \to x_0 + \eps v_2 \to x_0]] = 1_G + \frac{\eps^2}{2} F(A)(x_0)(v_1, v_2) + o(\eps^2);
\end{equation}
indeed, this can be taken to be a more fundamental definition of the curvature $F(A)$ than \eqref{f-def}.

It is easy to verify the inversion law
\be{inversion}
A[[x_0 \to x_1]] = A[[x_1 \to x_0]]^{-1}
\end{equation}
and the concatenation law
\be{concatenation}
A[[x_0 \to x_1 \to x_2]] = A[[x_0 \to x_2]]
\end{equation}
when $x_1$ lies in $[x_0 \to x_2]$; these allows us to perform manipulations such as
$$ A[[x_0 \to x_1 \to x_2 \to x_0]] A[[x_0 \to x_2 \to x_3 \to x_0]] = A[[ x_0 \to x_1 \to x_3 \to x_0]]$$
whenever $x_2$ lies in $[x_1 \to x_3]$.

By conjugating \eqref{A-def} by $\sigma$ it is easy to arrive at the gauge transformation law
\be{trans}
\sigma(A)[[x_0 \to x_1]] = \sigma(x_0) A[[x_0 \to x_1]] \sigma(x_1)^{-1}
\end{equation}
whenever $[x_0 \to x_1]$ is disjoint from $S$.  More generally we have
$$
\sigma(A)[[x_0 \to x_1 \to x_2]] = \sigma(x_0) A[[x_0 \to x_1 \to x_2]] \sigma(x_2)^{-1},$$
whenever $[x_0 \to x_1 \to x_2]$ is disjoint from $S$, etc.  

If $A$ has no curvature, $F(A) = 0$, then the monodromy along any loop is zero, and in particular
we have $A[[x_0 \to x_1 \to x_2 \to x_0]]=1_G$ whenever the triangle $\Delta(x_0,x_1,x_2)$ does not intersect $S$.
When the curvature is non-zero, we can still estimate the monodromy by the integral of the curvature:

\begin{lemma}[Nonabelian Stokes theorem]\label{stokes}  Let $x_0,x_1,x_2 \in [0,1]^n$ be such that the triangle 
$\Delta(x_0,x_1,x_2)$
is disjoint from $S$.  Then
\be{stokes-eq}
|A[[x_0 \to x_1 \to x_2 \to x_0]] - 1_G| \leq C \int_{\Delta(x_0,x_1,x_2)} |F(A)(y)|\ d\H^2(y)
\end{equation}
where $d\H^2(y)$ is two-dimensional Hausdorff measure.
\end{lemma}

\begin{proof}
It will be convenient to replace $|A[]-1_G|$ by a slightly different quantity.  Let $d(,)$ be the arclength metric on 
the group $G \subseteq U(N)$; this metric is bi-invariant under left and right-multiplication by elements of $G$, and $d(g,
g') \sim |g-g'|$.  Thus it will suffice to prove the estimate
$$ 
d(A[[x_0 \to x_1 \to x_2 \to x_0],1_G) \leq C\int_{\Delta(x_0,x_1,x_2)} |F(A)(y)|\ d\H^2(y).$$

The claim is invariant under cyclic permutation of $x_0$, $x_1$, $x_2$.  Also, if we let $x_3$ be the
midpoint of $x_0$ and $x_2$, then from the invariance properties of the metric and the triangle inequality we have
\bas 
d(A[[x_0 \to x_1 \to x_2 \to x_0],1_G) \leq &
d(A[[x_0 \to x_1 \to x_3 \to x_0],1_G) \\
&+ d(A[[x_3 \to x_2 \to x_1 \to x_3]],1_G)
\end{align*}
while we trivially have
\bas \int_{\Delta(x_0,x_1,x_2)} |F(A)(y)|\ d\H^2(y) = &\int_{\Delta(x_0, x_1, x_3)} |F(A)(y)|\ d\H^2(y)\\
&+ \int_{\Delta(x_3, x_2, x_1)} |F(A)(y)|\ d\H^2(y).\
\end{align*}
Thus to prove the claim for the triangle $\Delta(x_0,x_1,x_2)$ it suffices to do so for the two smaller
triangles $\Delta(x_0, x_1, x_3)$ and $\Delta(x_3, x_2, x_1)$ (with exactly the same constant $C$).  
Repeating this calculation
indefinitely we thus see that it will suffice to prove this estimate for infinitesimal triangles.  But this
follows directly from \eqref{f-tri}.
\end{proof}

Of course, all the above results only hold subject to the caveat that various lines and triangles do
not intersect $S$.  Fortunately, if one has enough free parameters then these intersections are extremely rare\footnote{A
variant of this argument can be used to show that $[0,1]^n \backslash S$ is simply connected; indeed, any closed loop in $[0,1]^n \backslash S$ can be contracted along a generic cone over that loop.  While we will not
use this fact directly, it does shed some light as to why gauge transform results such as Proposition \ref{inductive-gauge}
or Proposition \ref{m-approx} are possible.}:

\begin{lemma}\label{sing}  For any point $x_0 \not \in S$, the set 
$\{ x \in \R^n: [x_0 \to x] \hbox{ intersects } S \}$ has measure zero.  For any $x_0, x_1 \not \in S$ with
$[x_0 \to x_1]$ not intersecting $S$, the set $\{ x \in \R^n: \Delta(x_0,x_1,x) \hbox{ intersects } S\}$
has measure zero.
\end{lemma}

\begin{proof}
It suffices to prove the second claim, since the first follows by setting $x_0 = x_1$.

Fix $x_0, x_1$.  If $\Delta(x_0,x_1,x)$ intersects $S$, then we must have
$$ y = (1-\alpha) ((1-\theta) x_0 + \theta x_1) + \alpha x$$
for some $0 \leq \alpha, \theta \leq 1$ and $y \in S$.  Since $[x_0 \to x_1]$ is disjoint from the compact set $S$,
we see that $\alpha > c > 0$ for some $c = c(x_0, x_1, S)$.  We now solve for $x$ as
$$ x = \frac{1}{\alpha} y - \frac{1-\alpha}{\alpha} ((1-\theta) x_0 + \theta x_1).$$
Since $y$ lives in a set of dimension at most $n-4$, and $\alpha$ and $\theta$ are one-dimensional parameters
with $1/\alpha$ bounded, it is easy to see that $x$ lives in a set of dimension at most $n-2$, which necessarily
has measure zero, as desired.
\end{proof}

To exploit this generic lack of intersections, we shall use an averaging argument, using a random origin $x_0$ to
create a partially defined gauge (e.g. by using a radial gauge $\sigma(x) = A[[x_0 \to x]]$, defined as long as $[x_0 \to x]$ does not intersect $S$), and then averaging 
over $x_0$ (using Lemma \ref{sing}) to recover a globally defined gauge.

To do this we need a notion of averaging\footnote{One can view the machinery here as a continuous version of the more discrete gauge gluing techniques in, say, \cite{uhlenbeck}.  We were forced to use this continuous procedure instead of the discrete one in order to preserve the constants in the inductive procedure in Section \ref{inductive-sec}.} on the Lie group $G$.  Suppose that we have a domain $\Omega \subseteq [0,1]^n$,
a weight function $a: \Omega \to \R^+$  with $0 < \| a \|_{L^1(\Omega)} < \infty$, and a measurable map $f: \Omega \to G$ 
defined for a.e. $x \in \Omega$.  We would like to define a group element $[f]_{\Omega, a}^G \in G$ which represents in 
some sense an ``average'' of $f(x)$ where $x$ ranges over the probability measure $a dx/\|a\|_{L^1(\Omega)}$.

Since $G$ is embedded in the vector space $M_N(\C)$ of $N \times N$ complex matrices, we can define the linear average 
$[f]_{\Omega, a} \in M_N(\C)$ by
$$ [f]_{\Omega, a} := \frac{\int_\Omega f(x) a(x)\ dx}{\int_\Omega a(x)\ dx}$$
but of course this average will almost certainly lie outside of the group $G$.  To resolve this problem, we observe from the 
compactness of $G$ that there exists a tubular $\delta$-neighbourhood $N_\delta(G)$ of $G$ for some fixed 
$0 < \delta \ll 1$ for which there is a smooth projection map $\pi: N_\delta(G) \to G$ which equals the identity on 
$G$.  We can in fact choose $\pi$ to obey the right-equivariance condition $\pi(xg) = \pi(x)g$ for all 
$x \in N_\delta(G)$ and $g \in G$.  In particular we observe (since $\pi$ is Lipschitz) that we have the 
estimate
\be{lippi}
| \pi(x) g \pi(y)^{-1} - 1 | = | \pi(xg) \pi(y)^{-1} - 1 | \leq C |\pi(xg) - \pi(y) | \leq C |xg - y|
\end{equation}
for all $g \in G$ and $x,y \in N_\delta(G)$.

We now define
$$ [f]_{\Omega, a}^G := \pi( [f]_{\Omega, a} )$$
provided that the average $[f]_{\Omega, a}$ lies in $N_\delta(G)$.  This can be achieved provided that the values of 
$f$ ``cluster together''; more precisely, we have

\begin{lemma}\label{cluster-lemma} Let the notation and assumptions be as above.  If $f$ obeys the clustering condition
\be{cluster}
\int_\Omega \int_\Omega |f(x) - f(y)| a(x) a(y)\ d\mu(x) d\mu(y) <
\delta \| a \|_{L^1(\Omega)}^2
\end{equation}
then $[f]_{\Omega, a}$ lies in $N_\delta(G)$, and so $[f]^G_{\Omega, a}$ is well defined.  Here we use $d\mu$ to
denote Haar measure on the compact group $G$, normalized so that $\mu(G) = 1$.
\end{lemma}

\begin{proof}  We may of course assume that $f(x)$ is defined and takes values in $G$ for all $x \in X$, since sets of 
measure zero are clearly irrelevant here.
 
By \eqref{cluster} and the pigeonhole principle we may find an $x \in X$ such that 
$$
\int_\Omega |f(x) - f(y)| a(y)\ d\mu(y) <
\delta \| a \|_{L^1(\Omega)}.$$
From the triangle inequality we thus have
$$
|\int_\Omega (f(x) - f(y)) a(y)\ d\mu(y)| <
\delta \| a \|_{L^1(\Omega)}.$$
But the left-hand side simplifies to
$$ |f(x) - [f]_{\Omega, a}| \| a \|_{L^1(\Omega)}.$$
Thus $[f]_{\Omega, a}$ lies within distance $\delta$ of the point $f(x) \in G$, and the claim follows.
\end{proof}

\section{Fractional integration bounds}\label{frac-sec}

In the next section we shall be constructing gauges by averaging certain integrals of the form 
$A[[x_0 \to x_1 \to x_2]]$ using 
Lemma \ref{cluster-lemma}.  We will then use Lemma \ref{stokes} to estimate the expressions which then result.  
This will lead to integrating the curvature on an ``averaged collection of 2-surfaces''; to assist the proof of 
Proposition \ref{inductive-gauge} we now present a simple lemma (mostly a consequence of the change of variables 
formula) to understand such integrals.

\begin{lemma}\label{triangle-frac}  For any $0 < r \leq R \leq 1$ and $x \in [0,1]^n$, we have
\be{tf}
\begin{split}
\int_{B(x,R)} &\int_{B(x,r)} \int_{\Delta(x, x_1, x_2)} |F(A)(y)|\ d\H^2(y)\ dx_1 dx_2
\leq \\
&C r^n R^n (\int_{B(x,r)} |F(A)(y)| |x-y|^{2-n}\ dy \\
&+ r \int_{B(x,2R) \backslash B(x,r)} |F(A)(y)| |x-y|^{1-n}\ dy).
\end{split}
\end{equation}
\end{lemma}

\begin{proof}
We parameterize $x_1 = x + rz$, $x_2 = x + Rz'$, and $y = x + trz + t'Rz'$ for $z,z' \in B(0,1)$ and $t,t' \in [0,1]$ 
(we also have $t+t' \leq 1$, but we will not need this), and bound the left-hand side of \eqref{tf} by
$$ C r^n R^n r R \int_{B(0,1)} \int_{B(0,1)} \int_0^1 \int_0^1
|F(A)(x + trz + t'Rz')|\ dt dt' dz dz'.$$
Making the change of variables $(t,t',w,w') = (t,t',tz,tz')$, this becomes
$$ C r^n R^n r R \int_0^1 \int_0^1 \int_{B(0,t')} \int_{B(0,t)}
(tt')^n |F(A)(x + rw + Rw')|\ dw dw' dt dt'.$$
Applying Fubini's theorem and performing the $t$, $t'$ integrals, this becomes
$$ C r^n R^n r R \int_{B(0,1)} \int_{B(0,1)} |w|^{1-n} |w'|^{1-n} 
|F(A)(x + rw + Rw')|\ dw dw'.$$
Making the change of variables $(w,y) = (w, x+rw+Rw')$, and noting that $y \in B(x,2R)$, this can be bounded by
$$ C r^n R^n r R R^{-n} \int_{B(x,2R)} \int_{B(0,1)} (\frac{|y-x-rw|}{R})^{1-n} |w|^{1-n} 
|F(A)(y)|\ dw dy.$$
Comparing this with \eqref{tf}, it thus suffices to show that
$$ \int_{B(0,1)} |y-x-rw|^{1-n} |w|^{1-n}\ dw \leq C |y-x|^{2-n} / r$$
when $|y-x| \leq r$ and
$$ \int_{B(0,1)} |y-x-rw|^{1-n} |w|^{1-n}\ dw \leq C |y-x|^{1-n}$$
otherwise.  But this can be verified by a direct computation.
\end{proof}

From the above lemma we see that it is of interest to compute various fractional integrals of $F(A)$; we record two 
such computations below.

\begin{lemma}\label{fract-int}  If $x \in \Omega_m$ and $r > 0$, then
\be{om-1}
\int_{B(x,r)} |F(A)(y)| |x-y|^{2-n}\ dy \leq C \eps \log(2 + \frac{r}{R_m}).
\end{equation}
We also have the variant estimate: if $x \in \Omega_{m+1}$, then
\be{om-2}
\int_{[0,1]^n} 
\min(1, (\frac{R_m}{|x-y|})^{1-\kappa/2})
\frac{|F(A)(y)|}{|x-y|^{n-2}}\ dy \leq C \eps \log D.
\end{equation}
\end{lemma}

\begin{proof}
We decompose the left-hand side of \eqref{om-1} dyadically, and bound it by
\be{cs}
C \sum_{j=0}^\infty (2^{-j} r)^{2-n} \int_{B(x, 2^{-j} r)} |F(A)|.
\end{equation}
By \eqref{F-avg} we have the bound
$$
\int_{B(x, 2^{-j} r)} |F(A)| \leq C \eps (2^{-j} r)^{n-2}.
$$
On the other hand, since $x \in \Omega_m$, we have $Q(x) \leq \eps D^{(1-\kappa)m} = C \eps R_m^{-1+\kappa}$, which 
implies from \eqref{Q-def} and H\"older that
$$
\int_{B(x, 2^{-j} r)} |F(A)| \leq 
C (2^{-j} r)^{n/2} (\int_{B(x, 2^{-j} r)} |F(A)|^2)^{1/2}
\leq
C \eps (2^{-j} r)^{n-1-\kappa} R_m^{-1+\kappa}.
$$
Thus we can estimate the left-hand side of \eqref{cs} by
$$ C \sum_{j=0}^\infty \min(\eps, \eps (2^{-j} r / R_m)^{1-\kappa}),$$
and the claim \eqref{om-1} follows.

Now we prove \eqref{om-2}.  
The portion of the integral when $|x-y| \leq R_m$ is acceptable by \eqref{om-1} (with $m$ replaced by $m+1$,
and recalling that $R_m/R_{m+1} = D$).  
The portion when $|x-y| > R_m$ can be decomposed dyadically, and estimated by
$$ C\sum_{j = 0}^\infty 2^{-j(1-\kappa/2)} (2^j R_m)^{2-n} \int_{B(x, 2^j R_m)} |F(A)|.$$
By \eqref{F-avg} this can be bounded by
$$ C\sum_{j=0}^\infty 2^{-j(1-\kappa/2)} \eps$$
which is acceptable.
\end{proof}

From the above lemma we see that the integral
$$ \int_{[0,1]^n} \frac{|F(A)(y)|}{|x-y|^{n-2}}\ dy$$
might diverge logarithmically, like $\eps |\log\rho(x)|$, as one approaches the singular set $S$.  This is rather unfortunate; if this integral were uniformly bounded by $O(\eps)$, then one would not need the rather complicated inductive argument below, as the $m=1$ iteration of the gauge would already extend all the way down to $S$.  It is interesting that Price's monotonicity formula \cite{price} does give some additional control on the radial component $F(A)(y) \cdot \frac{x-y}{|x-y|}$ of the curvature, in particular obtaining an estimate of the form
$$ \int_{[0,1]^n} \frac{| F(A)(y) \cdot \frac{x-y}{|x-y|}|^2}{|x-y|^{n-4}}\ dy \leq C \eps.$$
Furthermore, in all the arguments in the next section it turns out that we only need the radial component of the curvature.  Unfortunately, this bound only seems able to improve the logarithmic divergence slightly, to $O(\eps \sqrt{|\log \rho(x)|})$, but cannot eliminate it entirely.  Thus we have been forced to perform this somewhat artificial and complicated inductive procedure in order to obtain a gauge which extends arbitrarily close to the singular set $S$ and which obeys manageable bounds on the connection $\sigma_m(A)$. 

\section{Proof of Proposition \ref{inductive-gauge}}\label{inductive-sec}

Armed with all the above machinery we can now prove Proposition \ref{inductive-gauge}.

We proceed by induction on $m$.  For inductive purposes we shall need to make certain constants explicit; specifically, 
we shall need a large constant $C_1$ depending on $D$, $n$, $G$.  The precise inductive claim is as follows:

\begin{proposition}\label{inductive-claim}  For each $m \geq 1$, there exists a gauge $\sigma_m$ with the following 
Lipschitz property: whenever $x_0, x_1 \in \Omega_m$ and $r > 0$ are such that $|x_0 - x_1| \leq r \leq 10 R_m$ we 
have\footnote{Of course, the integrand is only defined if $[x_0 \to x_2 \to x_1]$ avoids $S$, but this will turn out
to be the case for almost every $x_2$.  See Lemma \ref{sing}.}
\be{lip-0}
\begin{split}
r^{-n} \int_{B(x_0, r) \cap [0,1]^n}& |\sigma_m(x_0) A[[x_0 \to x_2 \to x_1]] \sigma_m(x_1)^{-1} - 1_G|\ dx_2 \\
&\leq C_1 r (T_m(x_0) + T_m(x_1)).
\end{split}
\end{equation}
\end{proposition}

Suppose for the moment that Proposition \ref{inductive-claim} held.
Let $x_0 \in \Omega_m$.  If we set $x_1 = x_0 + \epsilon v$ and $r := \epsilon$ for some unit vector $v$ and some 
$\epsilon$ small enough that $B(x_0, 2\epsilon) \subseteq \Omega_m$, then from \eqref{trans} we have
$$ |\sigma_m(x_0) A[[x_0 \to x_2 \to x_1]] \sigma_m(x_1)^{-1} - 1| = |\sigma_m(A)[[x_0 \to x_2 \to x_1]] - 1|$$
whenever $x_2 \in B(x_0,r) \cap [0,1]^n$, which by \eqref{A-def} implies that
$$ |\sigma_m(x_0) A[[x_0 \to x_2 \to x_1]] \sigma_m(x_1)^{-1} - 1| \geq
C^{-1} \epsilon |v \cdot \sigma_m(A)(x_0)| + o(\epsilon).$$
Combining this with \eqref{trans} and \eqref{lip-0}, we obtain
$$ \epsilon |v \cdot \sigma_m(A)(x_0)| + o(\epsilon) \leq C_1 \epsilon (T_m(x_0) + T_m(x_0 + \epsilon v)).$$
Dividing by $\epsilon$ and then taking limits as $\epsilon \to 0$, we obtain \eqref{sigma_m-bound-alt} as desired, since $v$ 
was arbitrary (note that $T_m$ is continuous on $[0,1]^n \backslash S$, and in particular on $\Omega_m$).    Thus to prove Proposition \ref{inductive-gauge} it 
will suffice to prove Proposition \ref{inductive-claim}.  As with all inductions, this is done in two steps.

{\bf Step 1. The base case $m=1$.}

We first construct the initial gauge $\sigma_1$.  We remark that it would be relatively easy to start the induction 
if we knew that $[0,1]^n \backslash S$ contained a large ball, but we are making no assumptions on $S$ other than the 
dimension assumption and so cannot assume this.  Besides, the arguments we will use here will also motivate the 
inductive step below.

We first pick a good choice of origin $x_* \in \Omega_1$.
From Fubini's theorem we have
$$ \int_{\Omega_1} \int_{[0,1]^n} \frac{|F(A)(x)|}{|x-x_*|^{n-2}}\ dx dx_*
\leq C \int_{[0,1]^n} |F(A)(x)|\ dx.$$
From \eqref{large-intersection} and the pigeonhole principle, there thus exists $x_* \in \Omega_1$ such that
\be{xs-bound}
\int_{[0,1]^n} \frac{|F(A)(x)|}{|x-x_*|^{n-2}}\ dx 
\leq C \int_{[0,1]^n} |F(A)(x)|\ dx.
\end{equation}
Fix this $x_*$. For each $x_1 \in [0,1]^n$ we define the group element $f_{x,0}(x_1) \in G$ by
\be{f0-def}
 f_{x,0}(x_1) := A[ [x_* \to x_1 \to x] ].
\end{equation}
This group element is undefined if $[x_* \to x_1 \to x]$ intersects $S$, but from Lemma \ref{sing} we see that the set of
$x_1$ for which that occurs has measure zero.  Thus $f_{x,0}: [0,1]^n \to G$ is defined at almost every 
point of $[0,1]^n$, and there will be no difficulty integrating this function in the $x_1$ variable.

\begin{lemma}[Clustering bound for $\sigma_1$]\label{f-cluster}  For all $x \in \Omega_1$, we have the clustering bound
\be{cl-1}
\int_{[0,1]^n} \int_{[0,1]^n} |f_{x,0}(x_1) - f_{x,0}(x_2)|\ dx_1 dx_2
\leq C \int_{[0,1]^n} |F(A)(y)| |x-y|^{2-n} \ dy.
\end{equation}
\end{lemma}

\begin{proof} From \eqref{f0-def}, \eqref{inversion} we see that $f_{x_0}(x_1) - f_{x_0}(x_2)$ is
conjugate to $A[[x_* \to x_1 \to x \to x_2 \to x_*]] - 1_G$, and in particular that
$$ |f_{x,0}(x_1) - f_{x,0}(x_2)| = |A[ [x_* \to x_1 \to x \to x_2 \to x_*] ] - 1_G|.$$
The loop on the right-hand side is the boundary of $\Delta(x_*, x_1, x_2) \cup \Delta(x, x_1, x_2)$, where 
$\Delta(x,x_1,x_2)$ denotes the triangle with vertices $x, x_1, x_2$.  For almost every $x_1$, we see from
Lemma \ref{sing} that the triangles $\Delta(x_*, x_1, x_2)$ and $\Delta(x, x_1, x_2)$ are 
disjoint from $S$ for almost every $x_2$.  We may thus apply Lemma \ref{stokes}, and bound the left-hand side of
\eqref{cl-1} by
$$C \int_{[0,1]^n} \int_{[0,1]^n} (\int_{\Delta(x,x_1,x_2) \cup \Delta(x_*,x_1,x_2)} |F(A)(y)| d\H^2(y))\ dx_1 dx_2.$$
By Lemma \ref{triangle-frac} we can bound this by
$$C \int_{[0,1]^n} |F(A)(y)| (|x-y|^{2-n} + |x_*-y|^{2-n})\ dy.$$
The claim then follows from \eqref{xs-bound}.
\end{proof}

From \eqref{cl-1} and Lemma \ref{fract-int} we have in particular that 
\be{cl-1a}
\int_{[0,1]^n} \int_{[0,1]^n} |f_{x,0}(x_1) - f_{x,0}(x_2)|\ dx_1 dx_2
\leq C \eps \log D 
\end{equation}
for $x \in \Omega_1$.

We now define $\sigma_1$ by averaging $f_{x,0}$ on the cube $[0,1]^n$:
$$ \sigma_1(x) := [f_{x,0}]_{[0,1]^n, 1}^G = \pi([f_{x,0}]_{[0,1]^n, 1});$$
from \eqref{cl-1a} and Lemma \ref{cluster-lemma} we see that $\sigma_1$ is well-defined on $\Omega_1$, if $\eps$ is 
sufficiently small depending on $D$.  
We now verify the condition \eqref{lip-0}; in other words, we show that whenever
 $x_0, x_1 \in \Omega_1$ and $|x_0 - x_1| \leq r \leq 10 R_1$, we have the bound
\be{pi-lip}
\begin{split}
r^{-n} \int_{B(x_0, r) \cap [0,1]^n} |\pi([f_{x_0,0}]_{[0,1]^n, 1}) A[[x_0 \to x_2 \to x_1]] &
\pi([f_{x_1,0}]_{[0,1]^n, 1})^{-1} - 1_G|\ dx_2 \\
\leq &C_1 r (T_1(x_0) + T_1(x_1)).
\end{split}
\end{equation}

Fix $x_0$, $x_1$, $r$.  By \eqref{lippi}, we can bound the left-hand side of \eqref{pi-lip} by
$$ C r^{-n} \int_{B(x_0, r) \cap [0,1]^n}  |[f_{x_0,0}]_{[0,1]^n, 1} A[[x_0 \to x_2 \to x_1]] - 
[f_{x_1,0}]_{[0,1]^n, 1}|;$$
by the triangle inequality, this is bounded by
$$ C r^{-n} \int_{B(x_0, r) \cap [0,1]^n} \int_{[0,1]^n} |f_{x_0,0}(x) A[[x_0 \to x_2 \to x_1]] - 
f_{x_1,0}(x)|\ dx dx_2.$$
From \eqref{f0-def}, \eqref{inversion} we see that $f_{x_0,0}(x) A[[x_0 \to x_2 \to x_1]] - 
f_{x_1,0}(x)$ is conjugate to $A[[x \to x_0 \to x_2 \to x_1 \to x]] - 1_G$, thus we can rewrite
the previous as
$$ C r^{-n} \int_{B(x_0, r) \cap [0,1]^n}
\int_{[0,1]^n} |A[[x \to x_0 \to x_2 \to x_1 \to x]] - 1_G|\ dx dx_2.$$
Applying Lemma \ref{stokes} we can bound this integral by
\be{r-block}
C r^{-n} \int_{B(x_0, r) \cap [0,1]^n} \int_{[0,1]^n}
(\int_{\Delta(x,x_2,x_0) \cup \Delta(x,x_2,x_1)} |F(A)(y)|\ d\H^2(y)) \ dx dx_2.
\end{equation}
First consider the integral on $\Delta(x, x_2, x_0)$.  From Lemma \ref{triangle-frac} (and using the crude 
estimate $|x-y|^{2-n} \leq r |x-y|^{1-n}$ when $|x-y| \leq r$) we can bound this by
$$ C r \int_{[0,1]^n} |x_0-y|^{n-1} |F(A)(y)|\ dy \leq C(D) r T_1(x_0)$$
as desired.  The contribution of the triangles $\Delta(x, x_2, x_1)$ is similar, but with $T_1(x_0)$ replaced 
by $T_1(x_1)$.  This proves \eqref{pi-lip}.

This completes the proof of the base case $m=1$.

{\bf Step 2.  Induct from $m$ to $m+1$.}

We now assume that Proposition \ref{inductive-claim} has already been proven for $m$, and now construct a gauge 
$\sigma_{m+1}$ on $\Omega_{m+1}$ with the desired properties.  This will basically be the same argument as Step 1, 
but rescaled by a factor of $D^{-m}$ and with some minor technical changes.  The key point here is that when moving 
from $\sigma_m$ to $\sigma_{m+1}$, the bound on the right-hand side of \eqref{lip-0} worsens by a factor of about $D$, 
which will allow us to close the argument if $D$ is sufficiently large.

Let $x \in \Omega_{m+1}$, and consider the ball $B_x := B(x, R_m)$.  
Let $\psi_x$ be the cutoff function $\psi_x(y) := \psi( (y-x)/R_m )$, where $\psi$ is a non-negative bump function 
adapted to $B(0,2)$ which equals one on $B(0,1)$; thus $\psi_x$ equals one on $B_x$.  In particular, from 
\eqref{large-intersection} we have
\be{psi-large}
\| \psi_x \|_{L^1(\Omega_m)} \geq |B_x \cap \Omega_m| \geq C^{-1} R_m^n.
\end{equation}
We define the function $f_{x,m}$ on $\Omega_m$ by
\be{fm-def}
f_{x,m}(x_1) := \sigma_m(x_1) A[[x_1 \to x]];
\end{equation}
observe that Lemma \ref{sing} ensures that $f_x$ is defined almost everywhere on $\Omega_m$.

One of the key observations we will need is that if we vary $x$ by $O(R_{m+1})$, then $\psi_x$ only varies by at
most $O(R_{m+1}/R_m) = O(1/D)$; this $1/D$ gain will be crucial in allowing us to close the induction.  This is the main reason why we need a large dyadic base $D$ instead of just using the standard powers of two.

The analogue of Lemma \ref{f-cluster} is

\begin{lemma}[Clustering bound for $\sigma_m$]\label{f-cluster-m}  
For all $x \in \Omega_{m+1}$, we have the clustering bound
\be{cl-2}
\begin{split}
\int_{B(x, 5R_m)} \int_{B(x, 5R_m)} &|f_{x,m}(x_1) - f_{x,m}(x_2)|\ dx_1 dx_2\\
&\leq C C_1 R_m^{2n}  \int_{[0,1]^n} 
\min(1, (\frac{R_m}{|x-y|})^{1-\kappa/2})
\frac{|F(A)(y)|}{|x-y|^{n-2}}\ dy.
\end{split}
\end{equation}
\end{lemma}

\begin{proof} Let $x_1, x_2 \in B(x, 5R_m)$.
By \eqref{fm-def}, \eqref{inversion} we have
$$ |f_{x,m}(x_1) - f_{x,m}(x_2)|
= |\sigma_m(x_1) A[[x_1 \to x \to x_2]] \sigma_m(x_2)^{-1} - 1_G|$$
since the expressions inside the absolute values are conjugate.
We insert a dummy variable $x_3$, ranging over $B(x_1, 20 R_m) \cap [0,1]^n$, and average to obtain
$$ |f_{x,m}(x_1) - f_{x,m}(x_2)|
\leq C
R_m^{-n} \int_{B(x_1, 20 R_m) \cap [0,1]^n}
 |\sigma_m(x_1) A[[x_1 \to x \to x_2]] \sigma_m(x_2)^{-1} - 1_G|\ dx_3.$$
We use the groupoid properties again and the triangle inequality to estimate the right-hand side by the sum of
\be{x-bypass}
C R_m^{-n} \int_{B(x_1, 20 R_m) \cap [0,1]^n}
|\sigma_m(x_1) A[[x_1 \to x_3 \to x_2]] \sigma_m(x_2)^{-1} - 1_G|\ dx_3
\end{equation}
and
\be{x-circuit}
C R_m^{-n} \int_{B(x_1, 20 R_m) \cap [0,1]^n}
|A[[x_1 \to x_3 \to x_2 \to x \to x_1]] - 1_G|\ dx_3.
\end{equation}

Consider first the contribution of \eqref{x-bypass} to \eqref{cl-2}.
Since $x_1, x_2 \in B(x, 5R_m)$, we have $|x_1 - x_2| 
\leq 10 R_m$.  Thus the inductive hypothesis \eqref{lip-0} applies, and we have the bound
$$
\eqref{x-bypass} \leq C C_1 R_m (T_m(x_1) + T_m(x_2)).$$
By symmetry, the contribution of \eqref{x-bypass} to the left-hand side of \eqref{cl-2} can thus be estimated by
$$ C C_1 R_m^n R_m \int_{B(x, 2R_m)} T_m(y)\ dy$$
which in turn can be estimated using \eqref{T_m-def} and Fubini's theorem by
$$ C C_1 R_m^n R_m \int_{[0,1]^n} \min(R_m , R_m^{n-\kappa/2} |x-y|^{1-n+\kappa/2}) |F(A)(y)|\ dy
$$
which is acceptable.

Now we consider the contribution of \eqref{x-circuit}.  By Lemma \ref{sing}, we see that the 
triangles $\Delta(x, x_1, x_3)$ and $\Delta(x_3, x_2, x)$ will not intersect $S$ for almost every choice of $x_1$, $x_2$, 
$x_3$.  Thus we may use Lemma \ref{stokes} and bound the contribution to \eqref{cl-2} by
\be{exp}
\begin{split}
C R_m^{-n} &
\int_{B(x, 5R_m)} \int_{B(x, 5R_m)} \int_{B(x_1, 20 R_m) \cap [0,1]^n} \\
&\int_{\Delta(x,x_1,x_3) \cup \Delta(x,x_2,x_3)} 
|F(A)(y)|\ d\H^2(y)
\ dx_3 dx_1 dx_2.
\end{split}
\end{equation}
Consider the $\Delta(x,x_1,x_3)$ integral.  By Lemma \ref{triangle-frac}, this portion of \eqref{exp} is bounded by
$$ C R_m^{2n} \int_{B(x, CR_m)} |x-y|^{2-n} F(A)(y)\ dy$$
which is acceptable.  The contribution of $\Delta(x,x_2,x_3)$ is similar.  This completes the proof of \eqref{cl-2}.
\end{proof}

From \eqref{cl-2}, \eqref{psi-large}, and \eqref{om-2} we have in particular that
\be{cl-1-m}
\int_{\Omega_m} \int_{\Omega_m} |f_{x,m}(x_1) - f_{x,m}(x_2)|\ \psi_x(x_1) dx_1 \psi_x(x_2) dx_2
\leq C C_1 \eps \| \psi_x \|_{L^1(\Omega_m)}^2 \log D.
\end{equation}
We now define the gauge $\sigma_{m+1}$ on $\Omega_{m+1}$ by the formula
$$ \sigma_{m+1}(x) := [f_{x,m}]_{\Omega_m, \psi_x}^G = \pi([f_{x,m}]_{\Omega_m, \psi_x}).$$
From \eqref{cl-1-m} and Lemma \ref{cluster-lemma} we see that
$\sigma_{m+1}$ is well-defined on $\Omega_{m+1}$.  Now we prove \eqref{lip-0}.

Fix $x_0, x_1 \in \Omega_{m+1}$ and suppose that $|x_0 - x_1| \leq r \leq 10 R_{m+1}$.  We have to prove that
\bas
r^{-n} &\int_{B(x_0, r) \cap [0,1]^n} |\pi([f_{x_0,m}]_{\Omega_m, \psi_{x_0}}) A[[x_0 \to x_2 \to x_1]] 
\pi([f_{x_1,m}]_{\Omega_m, \psi_{x_1}})^{-1} - 1|\ dx_2 \\
&\leq C_1 r (T_{m+1}(x_0) + T_{m+1}(x_1)).
\end{align*}
By \eqref{lippi}, we can bound the left-hand side by
$$
C r^{-n} \int_{B(x_0, r) \cap [0,1]^n} |[f_{x_0,m}]_{\Omega_m, \psi_{x_0}} A[[x_0 \to x_2 \to x_1]]  -
[f_{x_1,m}]_{\Omega_m, \psi_{x_1,m}}|\ dx_2.$$
By the triangle inequality it thus suffices to prove
\be{split-2}
\begin{split}
C r^{-n} \int_{B(x_0, r) \cap [0,1]^n} &| [f_{x_1,m}]_{\Omega_m, \psi_{x_0,m}} -  [f_{x_1,m}]_{\Omega_m, \psi_{x_1,m}}|\ 
dx_2 \\
&\leq \frac{1}{2} C_1 r (T_{m+1}(x_0) + T_{m+1}(x_1))
\end{split}
\end{equation}
and
\be{split-1}
\begin{split}
C r^{-n} \int_{B(x_0, r) \cap [0,1]^n} &|[f_{x_0,m}]_{\Omega_m, \psi_{x_0}} A[[x_0 \to x_2 \to x_1]]  -
[f_{x_1,m}]_{\Omega_m, \psi_{x_0,m}}|\ dx_2\\
&\leq \frac{1}{2} C_1 r (T_{m+1}(x_0) + T_{m+1}(x_1)).
\end{split}
\end{equation}
Let us first consider \eqref{split-2}.  We can integrate out the $x_2$ variable to bound the left-hand side by
$$
C | [f_{x_1,m}]_{\Omega_m, \psi_{x_0,m}} -  [f_{x_1,m}]_{\Omega_m, \psi_{x_1,m}}|$$
which we expand as
\be{tri}
C | \int_{\Omega_m} f_{x_1,m}(x)
\varphi(x)\ dx |
\end{equation}
where $\varphi$ is the function
$$ \varphi(x) := \frac{\psi_{x_0,m}(x)}{\int_{\Omega_m} \psi_{x_0,m}} - \frac{\psi_{x_1,m}(x)}{\int_{\Omega_m} 
\psi_{x_1,m}}.$$
From the support of $\varphi$ we may assume that $x$ is contained in the ball $B(x_0, 5R_m)$.

Since $\varphi(x)$ has mean zero, we may rewrite \eqref{tri} as
$$ C |\int_{\Omega_m} (f_{x_1,m}(x) - f_{x_1,m}(x'))
\varphi(x)\ dx |$$
for any $x'$. Averaging over $x' \in B(x_0, 5R_m)$ and using the triangle inequality, we may thus bound \eqref{tri} by
$$ C R_m^{-n} \int_{\Omega_m \cap B(x_0, 5R_m)} \int_{\Omega_m \cap B(x_0, 5R_m)}
| f_{x_1,m}(x) - f_{x_1,m}(x') | |\varphi(x)|\ dx dx'$$
which by Lemma \ref{f-cluster-m} is bounded by
$$ C C_1 R_m^n \| \varphi \|_\infty \int_{[0,1]^n} 
\min(1, (\frac{R_m}{|x_0-y|})^{1-\kappa/2})
\frac{|F(A)(y)|}{|x_0-y|^{n-2}}\ dy.$$
On the other hand, using the definition of $T_{m+1}(x_0)$  in \eqref{T_m-def} we have the bound
\bas
\int_{[0,1]^n} 
&\min(1, (\frac{R_m}{|x_0-y|})^{1-\kappa/2})
\frac{|F(A)(y)|}{|x_0-y|^{n-2}}\ dy\\
&\leq C D^{1-\kappa/2}
\int_{[0,1]^n} 
\min(1, (\frac{R_{m+1}}{|x_0-y|})^{1-\kappa/2})
\frac{|F(A)(y)|}{|x_0-y|^{n-2}}\ dy\\
&\leq C D^{1-\kappa/2} R_{m+1} T_{m+1}(x_0).
\end{align*}
Thus we can bound \eqref{tri} by
$$ C C_1 R_m^n \| \varphi \|_\infty 
D^{1-\kappa/2} R_{m+1} T_{m+1}(x_0).$$
We now compute $\| \varphi \|_\infty$. From the mean-value theorem we have
$$ \psi_{x_1}(x) = \psi_{x_0}(x) + O(\frac{r}{R_m})$$
and hence
$$ \int_{\Omega_m} \psi_{x_1}(y)\ dy = \int_{\Omega_m} \psi_{x_0}(y)\ dy + O(\frac{r}{R_m}) R_m^n.$$
Since $r \leq 10 R_{m+1} \ll R_m$, we thus see from \eqref{psi-large} that
$$ \| \varphi \|_\infty \leq C \frac{r}{R_m} R_m^{-n} = \frac{C}{D} \frac{r}{R_{m+1}} R_m^{-n}.$$
Thus we can bound \eqref{tri} by
$$ C C_1 r D^{-\kappa/2} T_{m+1(x_0)}$$
which will be acceptable if $D$ is large enough.

Now we consider \eqref{split-1}.  By the triangle inequality and \eqref{psi-large}, we may bound the left-hand side by
$$
C r^{-n} R_m^{-n}
\int_{B(x_0, r) \cap [0,1]^n} \int_{B(x_0, 5R_m) \cap [0,1]^n} |f_{x_0,m}(x) A[[x_0 \to x_2 \to x_1]]  -
f_{x_1,m}(x)|\ dx\ dx_2.$$
By \eqref{fm-def}, \eqref{inversion}, and some algebra, this can be bounded by
$$
C r^{-n} R_m^{-n}
\int_{B(x_0, r) \cap [0,1]^n} \int_{ B(x_0, 5R_m) \cap [0,1]^n} |A[[x \to x_0 \to x_2 \to x_1 \to x]]  - 1_G|\ dx\ dx_2.$$
Lemma \ref{sing} shows that for almost every $x$, $x_2$, the triangles $\Delta(x_0,x_2,x)$ and 
$\Delta(x_1, x_2, x)$ do not intersect $S$.  Thus we may use Lemma \ref{stokes} to estimate the previous by
\be{exp-2}
C r^{-n} R_m^{-n}
\int_{B(x_0, r) \cap [0,1]^n} \int_{ B(x_0, 5R_m) \cap [0,1]^n} \int_{\Delta(x_0,x_2,x) \cup \Delta(x_1,x_2,x)}
|F(A)(y)|\ d\H^2(y) dx\ dx_2.
\end{equation}
(compare with \eqref{r-block}).  Consider the $\Delta(x_0,x_2,x)$ integral.  By Lemma \ref{triangle-frac} (using 
the crude bound $|x-y|^{2-n} \leq |x-y|^{1-n} r$ when $|x-y| \leq r$) we can bound this portion of \eqref{exp-2} by
$$ C r \int_{B(x,CR_m)} |x_0-y|^{n-1} |F(A)(y)|\ dy.$$
By \eqref{T_m-def} we can bound this by $C(D) r T_{m+1}(x_0)$, which is acceptable if $C_1$ is large enough
depending on $D$.  The contribution of $\Delta(x_1,x_2,x)$ is similar but uses $T_{m+1}(x_1)$ instead of $T_{m+1}(x_0)$.  
This proves \eqref{split-2}, and closes the inductive step.  The proof of Proposition \ref{inductive-gauge} is now 
complete.
 
\section{Proof of Proposition \ref{m-approx}}\label{approx-sec}

In this section we use Proposition \ref{inductive-gauge} to prove Proposition \ref{m-approx}.

Fix $A$, $S$, $m$.  From Proposition \ref{inductive-gauge} we can find a gauge $\sigma_m$ on $\Omega_m$ obeying the bounds
\be{smax}
|\sigma_m(A)(x)| \leq C(D) \int_{[0,1]^n} (1 + \frac{|y-x|}{R_m})^{\kappa/2} 
\frac{|F(A)(y)|}{|x-y|^{n-1}}\ dy.
\end{equation}
for all $x \in \Omega_m$.

The gauge $\sigma_m(A)$ is currently only defined in $\Omega_m$.  The idea is now to truncate the gauge $\sigma_m(A)$ 
away from the set \eqref{open} to make it defined everywhere; then we use a mollification argument to make the truncated
gauge smooth.

We first need a Vitali covering argument.  For every $x \in [0,1]^n \backslash \Omega_m$, we see from \eqref{omega-def}
that there exists a radius $0 < r(x) \leq \diam([0,1]^n)$ such that
\be{ball-def}
r(x)^{-n/2+1+\kappa} (\int_{B(x,r(x))} |F(A)(y)|^2\ dy)^{1/2} \geq \frac{1}{2} \eps D^{(1-\kappa)m},
\end{equation}
and furthermore that
\be{ball-max}
r^{-n/2+1+\kappa} (\int_{B(x,r)} |F(A)(y)|^2\ dy)^{1/2} \leq 2 r(x)^{-n/2+1+\kappa} 
(\int_{B(x,r(x))} |F(A)(y)|^2\ dy)^{1/2}
\end{equation}
for all $0 < r \leq \diam([0,1]^n)$.
From \eqref{curv-small} and \eqref{ball-def} we see that 
\be{r-small}
r(x) \leq C D^{-m}
\end{equation}
for all $x \in [0,1]^n \backslash \Omega_m$.  

Since $[0,1]^n \backslash \Omega_m$ is compact, one can cover this space with only a finite number of balls $B(x, r(x))$.  
By the Vitali covering lemma, there thus exists a finite collection $x_1, \ldots, x_N$ of points in $[0,1]^n \backslash 
\Omega_m$ such that the balls $B(x_j, r(x_j))$ are disjoint, and that the balls $B(x_j, 5r(x_j))$ cover $[0,1]^n 
\backslash \Omega_m$.  In particular these balls also cover $S$.

For each $j = 1, \ldots, N$, let $\psi_j$ be a bump function adapted to $B(x_j, 10r(x_j))$ which equals one on $B(x_j, 
5r(x_j))$, and define the function  $\psi$ by
$$ \psi(x) := \sup_{j=1,\ldots,N} \psi_j(x).$$
Thus $\psi$ is a Lipschitz, piecewise smooth function which equals 1 on an open neighbourhood of $[0,1]^n \backslash
 \Omega_m$.  Since $r(x) \leq C D^{-m} \ll R_m$ we see from \eqref{omega-contains} that $\psi$ vanishes on the set 
$\{ x \in [0,1]^n: \rho(x) \geq 20 R_m \}$.  (Here of course we are using the fact that the connection $A$ is assumed
to obey the second property \eqref{f-point} of Proposition \ref{A-ass}, which is used to prove \eqref{omega-contains}).

We define the preliminary gauge $\tilde A_m$ on $[0,1]^n$ by the formula
\be{tam-def} 
\tilde A_m := (1-\psi) \sigma_m(A);
\end{equation}
note that even though $\sigma_m(A)$ is only defined on $\Omega_m$, $\tilde A_m$ is defined on all of $[0,1]^n$ since 
$1-\psi$ vanishes on a neighbourhood of $[0,1]^n \backslash \Omega_m$.  Also, $\tilde A_m$ is clearly gauge equivalent 
to $A$ on the set $\{ x \in [0,1]^n: \rho(x) \geq 20 R_m \}$, and in particular is smooth on this region.  Note also
that $\tilde A_m$ vanishes on a neighbourhood of $S$.

\begin{lemma}[$\tilde A_m$ has small curvature]  We have the curvature estimate
\be{fam-small}
\| F(\tilde A_m) \|_{M^{n/2}_2([0,1]^n)} \leq C(D) \eps.
\end{equation}
\end{lemma}

\begin{proof}
We first observe from \eqref{f-def}, \eqref{e-def} and the product rule that we have the pointwise estimate
\be{bound}
 |F(\tilde A_m)| \leq C |F(A)| + |\nabla \psi| |\sigma_m(A)|.
\end{equation}
Since $\nabla \psi$ is supported on $\Omega_m$ and obeys the pointwise bound $|\nabla \psi| \leq \sup_{j=1}^N |\nabla 
\psi_j|$, it thus suffices to show that
$$
\| \sup_j |\nabla \psi_j| |\sigma_m(A)| \|_{M^{n/2}_2(\Omega)} \leq C(D) \eps$$
or equivalently that
\be{bro}
\int_{B(x,r) \cap \Omega_m} \sup_j |\nabla \psi_j(x')|^2 |\sigma_m(A)(x')|^2\ dx' \leq C(D) \eps^2 r^{n-4}
\end{equation}
for all balls $B(x, r)$.

Fix $x$, $r$.  We may replace $\sup_j$ by $\sum_j$.  Since $\nabla \psi_j$ is supported on $B(x_j, 10r(x_j))$ and has 
magnitude $O(1/r(x_j))$, we can bound the left-hand side of \eqref{bro} by
$$ C \sum_j r(x_j)^{-2}
\int_{B(x,r) \cap B(x_j, 10r(x_j)) \cap \Omega_m}  |\sigma_m(A)(x')|^2\ dx';$$
applying \eqref{smax}, we thus reduce to showing that
\be{rsum}
\sum_j r(x_j)^{-2}
\int_{B(x,r) \cap B(x_j, 10r(x_j)) \cap \Omega_m}  
(\int_{[0,1]^n} (1 + \frac{|y-x'|}{R_m})^{\kappa/2} 
\frac{|F(A)(y)|}{|x'-y|^{n-1}}\ dy)^2 \ dx' \leq C(D) \eps^2 r^{n-4}.
\end{equation}
We now split into two cases, depending on whether $r(x_j) \geq r$ or $r(x_j) < r$.  First consider the terms where
$r(x_j) \geq r$.   From \eqref{r-small} this case can only occur when $r \leq C R_m$.  Now for each $k \geq 0$,
there are at most $O(1)$ balls $B(x_j, 10r(x_j))$ with $r(x_j) \sim 2^k r$ which intersect $B(x,r)$, since the
balls $B(x_j, r(x_j))$ are disjoint.  Thus we can sum the series $r(x_j)^{-2}$ and estimate this contribution
to the left-hand side of \eqref{rsum} by
$$ C r^{-2} \int_{B(x,r)}
(\int_{[0,1]^n} (1 + \frac{|y-x'|}{R_m})^{\kappa/2} 
\frac{|F(A)(y)|}{|x'-y|^{n-1}}\ dy)^2 \ dx'.$$
From \eqref{F-avg} and
a dyadic decomposition we have
$$ \int_{[0,1]^n: |x'-y| \geq r} (1 + \frac{|y-x'|}{R_m})^{\kappa/2} 
\frac{|F(A)(y)|}{|x'-y|^{n-1}}\ dy \leq C \eps/r$$
thus we may bound the previous by
$$ C \eps^2 r^{n-4} + C r^{-2} \int_{B(x,r)} (\int_{B(x',r)} \frac{|F(A)(y)|}{|x'-y|^{n-1}}\ dy)^2\ dx'.$$
The variable $y$ is now restricted to the ball $B(x,2r)$.  Since the kernel $\frac{1}{|x|^{n-1}}$ has an $L^1$ norm of
$O(r)$ on $B(0,r)$, we can use Young's inequality to estimate the previous by
$$ C \eps^2 r^{n-4} + C (\int_{B(x,2r)} |F(A)(y)|^2\ dy).$$
But this is acceptable by \eqref{curv-small}.

Now we check the contribution to \eqref{rsum} of the case where $r(x_j) < r$, which forces $B(x_j,r(x_j)) \subseteq
B(x, 20r)$.  We split the $y$ integration into
$|y-x_j| \leq 20r(x_j)$ and $|y-x_j| > 20r(x_j)$.  By \eqref{r-small}, the contribution when 
$|y-x_j| \leq 20r(x_j)$ is bounded by
$$ C \sum_{j: B(x_j, r(x_j)) \subseteq B(x,20r)} r(x_j)^{-2}
\int_{B(x_j, 10r(x_j))}  
(\int_{B(x',30r(x_j))} \frac{|F(A)(y)|}{|x'-y|^{n-1}}\ dy)^2 \ dx'.$$
Applying Young's inequality as before, we bound this by
$$ C \sum_{j: B(x_j, r(x_j)) \subseteq B(x,20r)} 
\int_{B(x_j,40r(x_j))} |F(A)(y)|^2\ dy.$$
Applying \eqref{ball-max}, we can bound this by
\be{ball-sum} C \sum_{j: B(x_j, r(x_j)) \subseteq B(x,20r)} 
\int_{B(x_j,r(x_j))} |F(A)(y)|^2\ dy.
\end{equation}
Since the balls $B(x_j, r(x_j))$ are disjoint, this is bounded by
$$ C \int_{B(x, 20r)} |F(A)(y)|^2\ dy$$
and the claim follows from \eqref{curv-small}.

Now we consider the contribution when $|y-x_j| > 20 r(x_j)$.  In this case the $x'$ variable is essentially irrelevant,
and we can estimate this contribution by
$$
C \sum_{j: B(x_j,r(x_j)) \subseteq B(x, 20r)} r(x_j)^{n-2}  
(\int_{|y-x_j| > 20 r(x_j)} (1 + \frac{|y-x_j|}{R_m})^{\kappa/2} 
\frac{|F(A)(y)|}{|x_j-y|^{n-1}}\ dy)^2;
$$
from \eqref{r-small} we may crudely bound this by
\begin{equation}\label{intermediate}
C \sum_{j: B(x_j,r(x_j)) \subseteq B(x, 20r)} r(x_j)^{n-2}  
(\int_{|y-x_j| > 20 r(x_j)} (\frac{|y-x_j|}{r(x_j)})^{\kappa/2} 
\frac{|F(A)(y)|}{|x_j-y|^{n-1}}\ dy)^2.
\end{equation}

From \eqref{ball-max} and H\"older we have
$$ \int_{|y-x_j| \leq r} |F(A)(y)| \leq C r^{n-1-\kappa} r(x_j)^{-n/2+1+\kappa}
(\int_{B(x_j,r(x_j))} |F(A)(y)|^2\ dy)^{1/2}$$
for any $0 < r \leq \diam([0,1]^n)$.  From dyadic decomposition we thus have
\bas
\int_{|y-x_j| > 20 r(x_j)} &(1 + \frac{|y-x_j|}{R_m})^{\kappa/2} 
\frac{|F(A)(y)|}{|x_j-y|^{n-1}}\ dy\\
\leq &\sum_{k \geq 0}
C 2^{k\kappa/2}
(2^k r(x_j))^{-(n-1)}
(2^k r(x_j))^{n-1-\kappa}\\
& r(x_j)^{-n/2+1+\kappa} 
(\int_{B(x_j,r(x_j))} |F(A)(y)|^2\ dy)^{1/2}\\
\leq &
C r(x_j)^{-n/2+1} 
(\int_{B(x_j,r(x_j))} |F(A)(y)|^2\ dy)^{1/2}
\end{align*}
since the $k$ summation is convergent.  In particular we have
$$ \eqref{intermediate} \leq C (\int_{B(x_j,r(x_j))} |F(A)(y)|^2\ dy)^{1/2}.$$
Thus we can bound this contribution by \eqref{ball-sum}, which is acceptable as before.
This proves \eqref{fam-small}.
\end{proof}

We are almost done, except that $\sigma_m$, and hence the $\tilde A_m$, is not necessarily smooth.  Fortunately
this can be easily resolved by regularizing $\sigma_m$.  

Fix $m$.  Since $\sigma_m$ is locally Lipschitz on $\Omega_m$, it lies in the Sobolev space $W^{1,p}$ on
the support of $1-\psi$ for any $n < p < \infty$.  In particular we can create a sequence $\sigma_{m,j}$ of 
smooth gauges which converge strongly in $W^{1,p}$ to $\sigma_m$ as $j \to \infty$ on the support of 
$1-\psi$; note that $W^{1,p}$ functions are H\"older continuous and so there is no difficulty keeping 
$\sigma_{m,j}$ on the Lie group $G$.

Define the gauges $\tilde A_{m,j}$ by
$$ \tilde A_{m,j} = (1-\psi) \sigma_{m,j}(A).$$
Then by construction $\tilde A_{m,j}$ is smooth, and vanishes near $S$.  Now we compare the curvatures of 
$\tilde A_{m,j}$
and $\tilde A_m$.  First we use \eqref{composition} to rewrite
$$ \tilde A_{m,j} = (1-\psi) \tilde \sigma_{m,j}(\sigma_m(A))$$
where $\tilde \sigma_{m,j} := \sigma_{m,j} \sigma_m^{-1}$.  From \eqref{f-def} and the product 
rule we have
$$ F(\tilde A_{m,j}) = (1-\psi) F(\tilde \sigma_{m,j}(\sigma_m A)) - 
\psi(1-\psi) \tilde \sigma_{m,j}(\sigma_m(A)) \wedge \tilde \sigma_{m,j}(\sigma_m(A))$$
and similarly
$$ F(\tilde A_m) = (1-\psi) F(\sigma_m(A)) - 
\psi(1-\psi) \sigma_m(A) \wedge \sigma_m(A).$$
On the other hand, from \eqref{gauge-identity} we have
$$ F(\tilde \sigma_{m,j}(\sigma_m(A))) = \tilde \sigma_{m,j}^{-1} F(\sigma_m(A)) \tilde \sigma_{m,j}$$
and hence
$$ F(\tilde A_{m,j}) - \tilde \sigma_{m,j}^{-1} F(\sigma_m(A)) \tilde \sigma_{m,j}
= \psi(1-\psi) (\tilde \sigma_{m,j}(\sigma_m(A)) \wedge \sigma_{m,j}(\sigma_m(A))
- \tilde \sigma_{m,j}^{-1} (\sigma_m(A) \wedge \sigma_m(A)) \tilde \sigma_{m,j}.$$
Since $\tilde \sigma_{m,j}$ converges to the identity in $W^{1,p}$ and $\sigma_m(A)$ is bounded on
the support of $1-\psi$, we see that the right-hand side converges to zero in $L^p$, and hence in $M^{n/2}_2$.
From \eqref{fam-small} we thus see that $\| F(\tilde A_{m,j})\|_{M^{n/2}_2([0,1]^n)} \leq \eps$ if
$j$ is sufficiently small.  The claim then follows by setting $A := \tilde A_{m,j}$ for this value of $j$.
The proof of Proposition \ref{m-approx} is now complete.

\section{Proof of Theorem \ref{abstract}}\label{setup-sec}

We now prove Theorem \ref{abstract}.  Our arguments are essentially those in 
\cite{uhlenbeck}, but with Lebesgue spaces replaced by their Morrey 
counterparts.  A similar argument has appeared in \cite{riviere}.

Let $K \gg 1$ be an absolute constant to be chosen later, and let $0 < \eps \ll 
1$ be sufficiently small depending on $K$.  Define $\U_\eps^*$ denote the space 
of connections $A \in \U_\eps$ which are gauge equivalent (via a smooth gauge) to a smooth Coulomb gauge 
$A_{coulomb}$ such that
\be{kb}
\| A_{coulomb} \|_{M^{n/2}_{2,1}([0,1]^n)} \leq K \eps.
\end{equation}
Clearly $\U_\eps^*$ is invariant under gauge transformations.

The main estimate is already contained in

\begin{lemma}[Bootstrap estimate]\label{main-lemma}
For any $A \in \U_\eps^*$, we can bootstrap \eqref{kb} to
\be{kb-2}
\| A_{coulomb} \|_{M^{n/2}_{2,1}([0,1]^n)} \leq K \eps/2.
\end{equation}
Also, we have \eqref{bound1}.
\end{lemma}

\begin{proof}
From elliptic estimates (Proposition \ref{elliptic-prop}) we 
have
$$ \| A_{coulomb} \|_{M^{n/2}_{2,1}([0,1]^n)} \leq C 
\|dA_{coulomb}\|_{M^{n/2}_2([0,1]^n)}.$$
By \eqref{f-def}, \eqref{hybrid-y}, \eqref{gauge-eq},  we thus have
$$ \| A_{coulomb} \|_{M^{n/2}_{2,1}([0,1]^n)}
\leq C (\|F(A)\|_{M^{n/2}_2([0,1]^n)} + \| A_{coulomb} 
\|_{M^{n/2}_{2,1}([0,1]^n)}^2).$$
The claims then follow from \eqref{kb}, \eqref{f-small} if $K$ is sufficiently 
large and $\eps$ sufficiently small depending on $K$.
\end{proof}

Fix $A \in \U_\eps$.  To prove Theorem \ref{abstract}, it suffices from the 
above lemma to show that $A \in \U_\eps^*$.

We exploit the smoothness of $A$ by choosing an exponent $n/2 < p < n$ (e.g. 
$p=3n/4$ will do).  We introduce the one-parameter family of connections $A_t$ 
for $t \in [0,1]$ by
$$ A_t(x) := t A(tx).$$
One may easily verify that the $A_t$ lie in $\U_\eps$, and that the map $t 
\mapsto A_t$ is continuous in the $M^p_{2,1}$ topology.  Also, $A_0 = 0$ is 
clearly in $\U^*_\eps$.  In order to show that $A_1$ is in $\U^*_\eps$ it thus 
suffices by standard continuity arguments to prove

\begin{proposition}[Continuity of the Coulomb gauge construction in smooth 
norms]\label{continuity}  Let $0 < X < \infty$, $n/2 < p < n$, and let $A \in 
\U^*_\eps$ be such that
\be{x-bound}
\| A \|_{M^p_{2,1}([0,1]^n)} \leq X.
\end{equation}
Then there exists a quantity $\delta_X > 0$ depending only on $X$, $G$, $n$, 
$p$, $\eps$, such that
$$
\{ A + \lambda \in \U_\eps: \| \lambda \|_{M^p_{2,1}([0,1]^n)} \leq \delta_X \}
\subset \U^*_\eps.
$$
\end{proposition}

\begin{proof}
Fix $p$; all our constants are allowed to depend on $p$.  We use $C_X$ to denote 
quantities which also depend on $X$.

{\bf Step 1.  Estimate the Coulomb gauge in smooth norms.}

From \eqref{x-bound}, \eqref{f-def}, \eqref{algebra} we have
$$ \| F(A) \|_{M^p_2([0,1]^n)} \leq C_X.$$
From \eqref{gauge-eq} we thus have
$$ \| F(A_{coulomb}) \|_{M^p_2([0,1]^n)} \leq C_X.$$
From elliptic estimates (Proposition \ref{elliptic-prop}) and 
\eqref{f-def} we thus have
$$ \| A_{coulomb} \|_{M^p_{2,1}([0,1]^n)} \leq C_X + C \| A_{coulomb} \wedge 
A_{coulomb} \|_{M^p_2([0,1]^n)}).$$
By \eqref{hybrid-x}, \eqref{kb} we thus obtain
\be{ac-bound}
\|A_{coulomb}\|_{M^p_{2,1}([0,1]^n)} \leq C_X.
\end{equation}

Next, let $\sigma_{coulomb}$ be the gauge transformation such that $A_{coulomb} 
= \sigma_{coulomb}(A)$.
From \eqref{e-def} we have
\be{e2}
d\sigma_{coulomb} = \sigma_{coulomb} A - A_{coulomb} \sigma_{coulomb}.
\end{equation}
Differentiating this, we have the pointwise estimate
$$ |\nabla^2 \sigma_{coulomb}| \leq C(
|\sigma_{coulomb}| (|\nabla A| + |\nabla A_{coulomb}|)
+ |d\sigma_{coulomb}|(|A|  + |A_{coulomb}|).$$
We substitute \eqref{e2} into this estimate.  Since $G$ is compact, 
$\sigma_{coulomb}$ is bounded,
and we thus obtain the pointwise bounds
$$ |\nabla^2 \sigma_{coulomb}|
\leq C( 1 + |\nabla A| + |A|^2 + |\nabla A_{coulomb}| +
|A_{coulomb}|^2).$$ From \eqref{e2} and the boundedness of
$\sigma_{coulomb}$ we in fact have
$$ |\sigma_{coulomb}| + |\nabla \sigma_{coulomb}| + |\nabla^2 \sigma_{coulomb}|
\leq C( 1 + |\nabla A| + |A|^2 + |\nabla A_{coulomb}| +
|A_{coulomb}|^2).$$ From \eqref{x-bound}, \eqref{ac-bound}, and
\eqref{algebra} we thus have \be{sc} \| \sigma_{coulomb}
\|_{M^p_{2,2}([0,1]^n)} \leq C_X.
\end{equation}

{\bf Step 2.  Pass to the Coulomb gauge.}

From \eqref{e-def}, \eqref{sc}, \eqref{algebra}, we see that the
gauge transformation $\tilde A \mapsto \sigma_{coulomb}(\tilde
A)$ is uniformly continuous in a small neighbourhood of $A$ (in
the $M^p_{2,1}$ topology).  Since $\U^*_\eps$ and $\U_\eps$ are
invariant under gauge transformations, we thus see that to prove
Proposition \ref{continuity} it suffices to do so in the case $A
= A_{coulomb}$.

{\bf Step 3.  Apply perturbation theory to the Coulomb gauge.}

Fix $\lambda$ as in the Proposition.  In order to show that $A+\lambda \in 
\U^*_\eps$, we first need to construct a gauge $\sigma$ such that
$$ d_*(\sigma(A + \lambda)) = 0 \hbox{ on } [0,1]^n$$
$$ n \cdot (\sigma(A + \lambda)) = 0 \hbox{ on } \partial [0,1]^n.$$
We use a perturbative argument.
Write $\sigma = \exp(U)$.  By \eqref{e-def} the above non-linear elliptic 
problem can then be written as
$$ \Delta U = d_* F(U, A+\lambda) \hbox{ on } \Omega$$
$$ n_\alpha d^\alpha U = F(U, A+\lambda) \hbox{ on } \partial \Omega$$
where
$$ F(U, A+\lambda) := dU - d\exp(U) \cdot \exp(-U) + \exp(U)\cdot (A + 
\lambda)\cdot \exp(-U).$$

We therefore set up the iteration scheme
$$ \Delta U^{(j+1)} = d_* F(U^{(j)}, A+\lambda) \hbox{ on } \Omega$$
$$ n_\alpha d^\alpha U^{(j+1)} = F(U^{(j)}, A+\lambda) \hbox{ on } \partial 
\Omega$$
$$ \int U^{(j+1)} = 0$$
where $U^0 := 0$.  Note that $U^{(j+1)}$ is uniquely defined by the classical 
theory of the Neumann problem; indeed, $U^{(j+1)}$ can be defined in terms of 
the Hodge decomposition of $F(U^{(j)}, A + \lambda)$.

We derive some bounds on $U^{(j+1)}$.  From the boundary conditions on $A$ we 
have $F(U^{(j)},\lambda) = F(U^{(j)},A+\lambda)$ on $\partial [0,1]^n$.  From 
Neumann problem regularity (Proposition \ref{neumann-prop}) we 
thus have
$$
\| U^{(j+1)} \|_{M^p_{2,2}([0,1]^n)} \leq C \| d_* F(U^{(j)}, A+\lambda) 
\|_{M^p_2([0,1]^n)}
+ C \| F(U^{(j)}, \lambda) \|_{M^p_{2,1}([0,1]^n)}.$$
From the chain rule and Definition \ref{coulomb-def} we observe the pointwise 
bounds
$$ |d_* F(U^{(j)}, A+\lambda)| \leq C( |U^{(j)}| |\nabla^2 U^{(j)}| + |\nabla 
U^{(j)}|^2 + |A| |\nabla U^{(j)}| + |\lambda| |\nabla U^{(j)}| + |\nabla 
\lambda| ).$$
By \eqref{algebra}, \eqref{bound1}, \eqref{hybrid-x} we thus have
$$ \| d_* F(U^{(j)}, A+\lambda)\|_{M^p_2([0,1]^n)} \leq C( \| U^{(j)} 
\|_{M^p_{2,2}([0,1]^n)}^2 + K\eps \| U^{(j)}\|_{M^p_{2,2}([0,1]^n)}
+ \| \lambda \|_{M^p_{2,1}([0,1]^n)}).$$
Similarly we have
$$ \| F(U^{(j)}, \lambda)\|_{M^p_{2,1}([0,1]^n)} \leq C( \| U^{(j)} 
\|_{M^p_{2,2}([0,1]^n)}^2 + \| \lambda \|_{M^p_{2,1}}([0,1]^n)).$$
Thus we have
$$ \| U^{(j+1)} \|_{M^p_{2,2}([0,1]^n)} \leq C( \| U^{(j)} 
\|_{M^p_{2,2}([0,1]^n)}^2 + K\eps \| U^{(j)}\|_{M^p_{2,2}([0,1]^n)}
+ \| \lambda \|_{M^p_{2,1}}([0,1]^n)).$$
If $\delta_X$ is sufficiently small, we thus obtain inductively
$$\|U^{(j)}\|_{M^p_{2,2}([0,1]^n)} \leq C \delta_X.$$
Adapting this scheme to differences, we thus see that $U^{(j)}$ converges in 
$M^p_{2,2}([0,1]^n)$ to a solution $U$ with
$$\|U\|_{M^p_{2,2}([0,1]^n)} \leq C \delta_X.$$
In particular, from Corollary \ref{sobolev} we see that $U$ has some H\"older continuity.  Standard elliptic regularity theory can then be used to bootstrap this regularity, eventually concluding that $U$ is smooth.  Exponentiating this (using 
\eqref{algebra}) we thus obtain a smooth Coulomb gauge $\sigma(A + \lambda)$ with
$$ \| \sigma - 1 \|_{M^p_{2,2}([0,1]^n)}, \| \sigma^{-1} - 
1\|_{M^p_{2,2}([0,1]^n)} \leq C \delta_X,$$
where $1$ is the identity element of $G$.
From \eqref{e-def}, \eqref{algebra}, \eqref{x-bound} we thus 
have
$$ \| \sigma(A + \lambda) - A \|_{M^p_{2,1}([0,1]^n)} \leq C_X \delta_X$$
and hence
$$ \| \sigma(A + \lambda) - A \|_{M^{n/2}_{2,1}([0,1]^n)} \leq C_X \delta_X.$$
If $\delta_X$ is sufficiently small depending on $X$, we thus see from 
\eqref{kb-2} that
$$ \| \sigma(A + \lambda) \|_{M^{n/2}_{2,1}([0,1]^n)} \leq K \eps$$
as desired.
\end{proof}

The proof of Theorem \ref{abstract} is now complete.  This completes all the
steps necessary to prove Theorem \ref{abstract3}.

\end{document}